%% file: main_RSP.tex
\documentclass[11pt, a4paper, oneside, reqno]{amsart}

\input{commands}

\renewcommand{\epsilon}{\varepsilon}
\newcommand{\blank}[1]{\hspace*{#1}\linebreak[0]}
\usepackage{dsfont,amssymb,amsmath, graphicx,fancyhdr,mdframed}
\usepackage{subcaption} 

\renewcommand{\epsilon}{\varepsilon}
\usepackage[utf8]{inputenc} 
\usepackage[T1]{fontenc}    
\usepackage{hyperref}       
\usepackage{url}            
\usepackage{booktabs}       
\usepackage{amsfonts}       
\usepackage{nicefrac}       
\usepackage{xcolor}         
\usepackage{pifont}
\usepackage{enumitem}
\usepackage{algorithm}
\usepackage[noend]{algpseudocode}
\usepackage{comment}
\newcommand{\edit}[1]{\textcolor{blue}{#1}}
\usepackage{adjustbox}

\newcommand{\R}{\mathbb{R}}

\usepackage[disable,textsize=tiny]{todonotes}
\usepackage[textsize=tiny]{todonotes}

\newenvironment{customproof}[1][\proofname]%
  {\par\noindent\textbf{#1. }\ignorespaces}%
  {\hfill\qed\par}
\input{style_PM}

\title[Locally Linear Convergence for Nonsmooth Convex Optimization]{Locally Linear Convergence for Nonsmooth Convex Optimization via Coupled \\ Smoothing and Momentum}

\author[R. {Rahimi Baghbadorani}]{Reza {Rahimi Baghbadorani}\textsuperscript{1}}
\author[S. Grammatico]{Sergio Grammatico\textsuperscript{1}}
\author[P. {Mohajerin Esfahani}]{Peyman {Mohajerin Esfahani}\textsuperscript{1,2}}

\thanks{{The authors are with (1) Delft University of Technology and (2) the University of Toronto. This work was supported by the ERC grant TRUST-949796 and the NSERC Discovery grant RGPIN-2025-06544.}}

\begin{document}

\maketitle

\begin{abstract}   
We propose an adaptive accelerated smoothing technique for a nonsmooth convex optimization problem where the smoothing update rule is coupled with the momentum parameter. We also extend the setting to the case where the objective function is the sum of two nonsmooth functions. With regard to convergence rate, we provide the global (optimal) sublinear convergence guarantees of $\mathcal{O}(1/k)$, which is known to be provably optimal for the studied class of functions, along with a local linear rate if the nonsmooth term fulfills a so-call {\em locally strong convexity condition}. We validate the performance of our algorithm on several problem classes, including regression with the $\ell_1$-norm (the Lasso problem), sparse semidefinite programming (the \textit{MaxCut} problem), Nuclear norm minimization with application in model free fault diagnosis, and $\ell_1$-regularized model predictive control to showcase the benefits of the coupling. An interesting observation is that although our global convergence result guarantees~$\mathcal{O}(1/k)$ convergence, we consistently observe a practical transient convergence rate of $\mathcal{O}(1/k^2)$, followed by asymptotic linear convergence as anticipated by the theoretical result. This two-phase behavior can also be explained in view of the proposed smoothing rule. 

\end{abstract}

\section{Introduction}\label{introduction}
Objective functions with multiple nonsmooth terms are widely used in optimization-based control, system identification and machine learning. The broad applicability of these formulations across diverse domains calls for efficient algorithms capable of handling nonsmooth optimization problems.\\
For instance, the authors in \cite{annergren2012admm,pakazad2013sparse} consider model predictive control (MPC) with nonsmooth $\ell_1$ regularizers to obtain sparse control inputs. The authors in \cite{noom2024proximal} use combinations of $\ell_1$ and Nuclear norms to model fault detection in model-free dynamical systems. \cite{fu2006efficient} employ $\ell_2$–$\ell_1$ and $\ell_1$–$\ell_1$ formulations for image restoration; \cite{chambolle2011first} utilize the TV–$\ell_1$ model for image denoising; \cite{wang2013dictionary} apply the $\ell_1$–$\ell_1$ formulation to dictionary learning; and \cite{yang2011alternating} use Nuclear norms for low-rank matrix decomposition in graph neural networks with an $\ell_1$ regularizer.

{\textbf{Motivating example.}} Consider the MPC formulation for a class of constrained linear systems with uncertain state-delays \cite{hu2004model}:  
\begin{align}\label{mpc-delayed}
      x(k+1) &= \mathbf{A} x(k) + \mathbf{A_d} x(k - N_d(x_k)) + \mathbf{B} u(k), \nonumber\\
      N_d(k) &\in [1, \bar{N}_d]; \quad x(i) = x_i, \quad \forall i \in \{-\bar{N}_d, \dots, 0\}, \nonumber\\
      \|u_i\| &\le \bar{u}_i, \quad i \in \{ 1, \dots, m\},
\end{align}
where $k$ is the discrete-time index, $x \in \mathbb{R}^n$ represents the system state, $u \in \mathbb{R}^m$ is the input vector, $\bar{u}_i$ are the input constraints, $N_d(k)$ is the uncertain time delay, possibly varying with time, $\bar{N}_d$ is its upper bound, and the matrices $\mathbf{A}$, $\mathbf{A_d}$, and $\mathbf{B}$ define the system dynamics. \\
Based on an artificial Lyapunov function, a stabilizing condition depending on the upper bound of the uncertain state-delays, together with MPC, is presented in the form of a linear matrix inequality as \cite[Eq. (15)]{hu2004model}:
\begin{align} \label{delayed-nonsmooth-sdp-mpc}
& \min_{ {u}(k|k), \dots, {u}(k+N-1|k) }  \ J(k) \\
&\text{s.t.}  \ \text{system dynamics and constraints \eqref{mpc-delayed},} \notag 
\end{align}
\begin{align*}
& \blank{.5cm}
\begin{bmatrix}
{P} - \bar{{N}}_d {Q}_d & 0 & (\mathbf{A} + \mathbf{B}{K})^\top \\
0 & {Q}_d & \mathbf{A}_d^\top \\
\mathbf{A} + \mathbf{B}{K} & \mathbf{A}_d & {P}^{-1}
\end{bmatrix} \succ 0, \notag \\
& \blank{.5cm}
\begin{bmatrix}
{Y} & {K} \\
{K}^\top & {P} \mu^{-1}
\end{bmatrix} \succeq 0, 
\quad
{Y}_{ij} \le \bar{u}_j^2, \quad \forall i,j \in \{ 1, \dots, m\}, \notag
\end{align*}
where, $J$ is a nonsmooth objective function that can include an $\ell_1$ regularizer, e.g., $\sum_i \| \Delta u_i \|_1$, to enforce sparsity in input changes, allowing the controller to respond promptly to disturbances and system variations \cite{pakazad2013sparse,annergren2012admm}. We emphasize that, the MPC problem in \eqref{delayed-nonsmooth-sdp-mpc}, is complex due to the nonsmooth objective function and SDP constraints, and must be solved at every sampling time. 

This motivates fast optimization algorithms for nonsmooth convex problem of the following class:
\begin{equation}\label{main problem}
    \min_{x \in \mathbb{R}^n} F(x) = \min_{x \in \mathbb{R}^n} f(x) + h(x), 
\end{equation}
where the functions~$f$ and $h$ are proper, closed, convex, possibly nonsmooth, but prox-friendly (to be made precise later in the notation section, see \eqref{prox-operator}). 
A particular subclass of such problems is composite minimization in which one of these functions is smooth~\cite{beck2009fast,nesterov2013gradient,baghbadorani2024adaptive}. While the composite minimization problem has been extensively studied, the possibility of having multiple nonsmooth terms remains relatively unexplored. This is primarily due to the fact that the summation of two prox-friendly functions is not necessarily prox-friendly~\cite{pustelnik2017proximity}. Let us emphasize that the most common algorithms, such as subgradient \cite{nemirovskij1983problem}, mirror descent \cite{beck2003mirror}, and bundle methods \cite{schramm1992version} often suffer from a slow convergence rate of~$\mathcal{O}\left(1/\epsilon ^2\right)$ where $\epsilon$ is an apriori desired precision.

Following the seminal works by Nesterov~\cite{nesterov2005smooth,nesterov2007smoothing}, one can exploit the structure of the nonsmooth objective function to propose a smooth approximation and then deploy a first-order accelerated method to find the optimizer of the smooth approximation~\cite{nesterov1983method}. Remarkably, Nesterov's algorithm requires only $\mathcal{O}(1/\epsilon)$~iterations to reach an $\epsilon$ optimal solution. Inspired by this approach, several works have tried to improve and enhance the efficacy of smoothing algorithms in terms of complexity or to customize it for specific applications at hand~\cite{nesta2009fast,nesterov2005excessive,chambolle2011first,d2004direct}. In \cite{beck2012smoothing}, comparisons are made about the advantages and disadvantages of smoothing techniques and proximal-type methods. Based on \cite{nesterov2007smoothing}, the authors in \cite{d2014stochastic} utilize a stochastic smoothing technique to improve the scalability of the algorithm for semidefinite programming problems. The smoothing techniques in~\cite{yurtsever2019conditional} improve the convergence speed in comparison with conditional gradient algorithms. The authors in \cite{orabona2012prisma} study the smoothing technique for minimizing the sum of three functions, where two of which are nonsmooth. The paper~\cite{boct2015variable} studies an adaptive smoothing algorithm with a convergence rate of $\mathcal{O}(1/\epsilon\log(1/\epsilon))$. 

The closest existing work to our study is the adaptive smoothing technique in~\cite{tran2017adaptive} that enjoys~$\mathcal{O}(1/\epsilon)$ complexity, which matches the worst-case bound offered by Nesterov's algorithm. Table~\ref{table1} summarizes the convergence results of this study and compares them with those in the literature. Following the footsteps of Nesterov's smoothing technique, we propose a novel adaptive-smoothing technique where the smoothing parameter decreases with each iteration but as a function of the momentum, which retains the same global convergence rate while improving the asymptotic performance. 
\begin{table*}[t]
    \centering
    \caption{Comparison of algorithms for minimizing nonsmooth convex functions.}
    \tiny
    \begin{tabular}{cccccc}\hline
    {Algorithm} & {Smoothing technique} & prox-friendly assumption & {Local linear convergence rate} & {Global convergence rate} \\
    \hline & \\[-1.0em]
    Subgradient descent~\cite{nemirovskij1983problem} & \ding{55} & \ding{55}  & \ding{55} & $\mathcal{O}(1/\varepsilon^2)$\\
    Stochastic smoothing~\cite{d2014stochastic} & \ding{55} & \ding{55} & \ding{55} & $\mathcal{O}(1/\varepsilon^2)$\\
     Chambolle-Poc~\cite{chambolle2011first} & \ding{55} & \ding{51} & \ding{55} & $\mathcal{O}(1/\varepsilon)$\\
    Nesterov's smoothing~\cite{nesterov2005smooth,nesterov2007smoothing} & \ding{51} & \ding{51} & \ding{55} & $\mathcal{O}(1/\varepsilon)$\\
    Variable smoothing~\cite{boct2015variable} & \ding{51} & \ding{51} & \ding{55} & $\mathcal{O}(\log(1/\varepsilon)/\varepsilon)$\\
    Adaptive smoothing~\cite{tran2017adaptive} & \ding{51} & \ding{51} & \ding{55} & $\mathcal{O}(1/\varepsilon)$\\
    \\[-1.0em]
    \hline
    \\[-1.0em]
    Proposed adaptive smoothing \\ (Theorems \ref{thm:global convegence}--\ref{thm:linear}) & \ding{51} & \ding{51} & \ding{51} (under some assumptions) & $\mathcal{O}(1/\varepsilon)$\\
    \\[-1.0em]
    \hline
    \end{tabular}
    \label{table1}
\end{table*} 

{\bf Contributions:} The contributions of our work are summarized as follows:
\begin{enumerate}[label = (\roman*), leftmargin=6mm, itemsep=1mm, topsep = 1mm]
     \item \textbf{Global optimal sublinear convergence:} 
     We introduce an algorithm with an adaptive smoothing parameter coupled with the momentum term (Algorithm~\ref{alg:Algorithm1}). When the smoothing rule is modified to stay away from zero, we provide a global sublinear convergence rate of~$\mathcal{O}\left(1/\epsilon\right)$ that matches the optimal worst-case complexity bound for optimizing this class of nonsmooth functions~(Theorem~\ref{thm:global convegence}).

    \item \textbf{Locally linear convergence:} 
     When the nonsmooth term~$f(x)$ meets a so-called {$\infty$-locally strong convexity condition}, we show that Algorithm~\ref{alg:Algorithm1} enjoys a local linear convergence rate~(Theorem~\ref{thm:linear}). Combined with the global convergence result from the previous section, this implies that an appropriate initial condition for Algorithm~\ref{alg:Algorithm1} ensures a transient optimal sublinear convergence rate followed by an asymptotic linear convergence rate. 
     
     \item \textbf{Multiple nonsmooth terms:}     
     The proposed algorithm allows for multiple nonsmooth terms. Such settings can be computationally challenging as the ``prox-friendly'' property, a key feature in nonsmooth optimization, is not necessarily additive. Important applications falling into this category include model-free fault diagnosis, nonsmooth model predictive control, sparse regression and sparse semidefinite programming, which are the examples investigated in our numerical section to validate the performance of our proposed algorithm.
\end{enumerate}

To validate the theoretical results of this study, we implement and compare the performance of the proposed algorithm with different existing methods in the literature on various classes of problems including regression with the $\ell_1$-norm (the Lasso problem), sparse semidefinite programming (the \textit{MaxCut} problem), and Nuclear norm minimization. Regarding the above contributions, an observation consistently confirmed by these numerical results deserves attention.   

\textbf{An unexpected observation}: 
The smoothness parameter has a direct impact on the stepsize in accelerated algorithms, i.e., the smoother the function, the larger the stepsize. A general rule of thumb is that larger stepsizes lead to faster convergence, and with that in mind, we would expect to see a slower convergence rate for coupled smoothness parameters. It, however, turns out that this is not the case and the proposed adaptive smoothness parameter has a faster convergence rate of $\mathcal{O}(1/k^2)$ for its transient and asymptotic convergence of $\mathcal{O}(e^{-k})$, as opposed to the existing rate of $\mathcal{O}(1/k)$~\cite{tran2017adaptive}. This two-phase behavior with these rates is also evident in the proposed smoothing rule~(cf. Lemma~\ref{lem:rate}).


\textbf{Roadmap}. The paper is organized as follows: Section~\ref{Literature review} reviews Nesterov's smoothing technique for solving nonsmooth minimization. We propose our adaptive-smoothing algorithm in Section~\ref{section3} and provide the technical proof of theorems in Section~\ref{thechnicalproof}. Section~\ref{simulation} benchmarks our algorithm in several applications.

\textbf{Notation}. $\mathbb{R}^n$ shows the real vector space, we denote the standard inner product by $\langle \cdot, \cdot \rangle$ and $\ell_p$-norm by $\|\cdot\|_p$ (and by $\|\cdot\|$, we mean the Euclidean standard 2-norm). If $f$ is differentiable, $\nabla f(x)$ represents the gradient of $f$ at $x$. The function~$f$ is called $L$-Lipschitz for some $L>0$ if
\begin{equation*}
    |f(x) - f(y)| \leq L\|x-y\|, \qquad \forall x,y \in \mathbb{R}^n.
\end{equation*}
The function~$f$ is $\mu$-smooth if its gradient is $\mu$-Lipschitz, i.e., 
$\|\nabla f(x) - \nabla f(y)\| \leq \mu\|x-y\|$. 
The convex function~$f$ is called {$\rho$-locally strongly convex at $x^*$, if there exists $\varepsilon > 0$ so that the function~$f(x) - \frac{\rho}{2}\|x\|^2$ is convex over the ball $\mathbb{B}_\varepsilon (x^*)$. We call $f$ $\infty$-locally strongly convex if $f$ is $\rho$-locally strongly convex for any $\rho > 0$}. The Fenchel conjugate of $f$ is defined as 
    $f^*(x) := \sup\limits_{y} \langle x,y \rangle - f(y)$.
The \text{prox} operator for function $h$ is defined as
\begin{equation}\label{prox-operator}
    \text{prox}_{h}(x) := \arg \min_{u} h(u) + \dfrac{1}{2} \|u-x\|^2.
\end{equation}
A function is ``prox-friendly'' if the operator~\eqref{prox-operator} is available (computationally or explicitly). We also denote the gradient mapping of two convex functions~$f$ and $h$ by 
\begin{align}
\label{grad_mapping}
G^{f}_{\zeta h}(x) := \dfrac{1}{\zeta}\Big(x - \mathrm{prox}_{\zeta  h}\big(x - \zeta \nabla f(x)\big)\Big),   
\end{align}
where $\zeta$ is a positive scalar and has a stepsize interpretation. The gradient mapping is available if $f$ is differentiable and $h$ is prox-friendly.
\section{State of the Art on Nonsmooth Optimization}\label{Literature review}

In this section, we review the current state of the art in smoothing nonsmooth functions and discuss a possible challenge that may emerge in dealing with multiple smoothing terms. 

\subsection{Nesterov's smoothing technique}\label{smoothing technique}
Consider a possibly nonsmooth convex function $f$ that we assume to be prox-friendly and $L_f$-Lipschitz continuous. A key object in smoothing techniques is the Moreau envelope defined~as
\begin{equation}
    f_{\mu}(x) := \min_{y} f(y) + \dfrac{1}{2\mu} \|y-x\|^2.
    \label{smooth appr}
\end{equation}
The Moreau envelope is a smooth lower approximation of a function at every point, i.e., $f_{\mu}(x) \leq f(x)$ for all $x \in \R^n$. By definition, the objective function of the Moreau envelop~\eqref{smooth appr} and the prox operator~\eqref{prox-operator} are closely related, namely, the optimizer of \eqref{smooth appr} is $\text{prox}_{\mu f}(x)$. The following lemma indicates several known properties of the smooth approximation~\eqref{smooth appr} that are central for algorithms in this context. For brevity, we skip the proof and refer interested readers to~\cite{boct2015variable} for further details.  

\begin{lemma}[Smoothness regularity]\label{lem:Moreau}
    Let $f_\mu$ be defined as in~\eqref{smooth appr}. Then, the following holds:
\begin{enumerate}[label=(\roman*), itemsep = 0mm, topsep = 0mm, leftmargin = 7mm]
    \item  \label{lem-dual}
    {\bf Dual reformulation:} 
    $f_\mu (x) = \max\limits_z \left\{\langle x,z \rangle - f^*(z) - \dfrac{\mu}{2}\|z\|^2\right\}$.
 
    \item \label{lem-bound}
    {\bf Uniform bound:}\label{(i)}  $f_{\mu}(x) \leq f(x) \leq f_{\mu}(x) + \dfrac{\mu}{2}L_{f}^2$, where $L_f$ is the Lipschitz constant of $f$.

   \item \label{lem-grad}
   {\bf Gradient evaluation:} 
   \label{(ii)} $\nabla f_\mu(x) = \dfrac{1}{\mu}\left(x - \text{prox}_{\mu f}(x)\right).$

   \item \label{lem-smooth} 
   {\bf Smoothness:} $f_\mu (\cdot)$ is $\dfrac{1}{\mu}$-smooth, i.e., $\|\nabla f_\mu(x) - \nabla f_\mu(y)\| \le {1 \over \mu}\|x - y\|$. 
\end{enumerate}
\end{lemma}

Lemma~\ref{lem:Moreau} paves the way to a power smoothing technique as follows: The uniform bound~\ref{lem-bound} allows us to choose a minimum value for the smoothing parameter so that the smoothed function~$f_\mu$ remains in a desired $\epsilon$-vicinity of the original function~$f$. Thanks to the smoothness result in~\ref{lem-smooth} and under the assumption that $f$ is prox-friendly, one can apply Nesterov's accelerated algorithm using the gradient evaluation~\ref{lem-grad} and optimize~$f_\mu$. The choice of Nesterov's algorithm (Algorithm~\ref{NesAlg}) is justified by the fact that it is the fastest general convex optimization algorithm for smooth functions~\cite{nesterov1983method}. 


\begin{algorithm}
\renewcommand{\thealgorithm}{0}
\caption{Nesterov's accelerated method for $1/\mu$-smooth functions \cite{nesterov1983method}}
\label{NesAlg}
\begin{algorithmic}
\State Input: given initial conditions $y_1 = x_1$, $\beta_0 > 0$, and smoothing constant $\mu$. 
\State $y_{k+1} = x_k - \zeta \nabla f_{\mu}(x_k)$, \quad  with the stepsize $\zeta = \mu$ \vspace{1mm} 
\State $x_{k+1} = (1-\gamma_k)y_{k+1} + \gamma_k y_k$
\State $\beta_{k+1} = {1 \over 2}\big(1+\sqrt{1+4\beta_{k}^2}\big),  \gamma_k = \dfrac{1-\beta_k}{\beta_{k+1}}$ 
\end{algorithmic}
\end{algorithm}

This idea is first proposed in \cite{nesterov2005smooth} in which the smoothness parameter is proposed to be the constant $\mu = 2\epsilon/L_f^2$ where $\epsilon$ is an apriori desired precision. This yields an overall complexity of $\mathcal{O}(1/\epsilon)$ in terms of the precision parameter~$\epsilon$. More recently, \cite{tran2017adaptive} proposes an adaptive version of the smoothing technique. Our proposed adaptive smoothing technique also follows a similar spirit as in~\cite{tran2017adaptive}, but the key feature is to pair the update rule with the momentum parameter of the accelerated method. Before proceeding with that, we also wish to briefly comment on our motivation for another feature of our proposed algorithm that allows for two nonsmooth terms. 

Several studies have been devoted to developing methods for computing the prox-operator~\eqref{prox-operator} of a sum of two (multiple) functions, see the recent work~\cite{pustelnik2017proximity,abboud2017dual,chaux2009nested,yu2013decomposing} and the references therein. These methods typically provide various assumptions under which we have 
\begin{equation}\label{sum of prox func}
    \text{prox}_{f+h} = \text{prox}_f \circ \text{prox}_h, 
\end{equation}
where ``$\circ$'' denotes the mapping composition. However, these conditions are still restrictive, and the prox-operator of the sum of two prox-friendly functions may not have a closed-form solution and may be in general computationally demanding (e.g., sparse regression and semi-definite programming). Motivated by this, we develop our algorithm so that it allows for two (multiple) nonsmooth terms, i.e., in~\eqref{main problem} both functions $f$ and $h$ are possibly nonsmooth but prox-friendly. 
\section{Proposed Algorithm and Convergence Analysis}\label{section3}

This section presents the main algorithm of this paper and its global and local convergence results. 

\subsection{Algorithm}\label{algorithms}
The proposed method is given in Algorithm~\ref{alg:Algorithm1}. The algorithm follows Nesterov's accelerated method in Algorithm~\ref{NesAlg} with the difference in the choice of stepsize~$\zeta$, which is coupled with the momentum parameter at two consecutive steps~$\beta_k$ and $\beta_{k+1}$ (Lemma \ref{prop:Alg1_bound}). This connection is established by looking at the increment of the function as described in one of the preliminary lemmas in Section \ref{thechnicalproof} (Lemma~\ref{lem:increment}). 


\begin{algorithm}
\renewcommand{\thealgorithm}{1}
\caption{Adaptive accelerated smoothing method}
\label{alg:Algorithm1}
\begin{algorithmic}[]
\State Input: given initial conditions $y_1 = x_1$, $\beta_0 > 0$, and the to-be-defined sequence $\big(\mu_k\big)_{k\geq0}$.  \vspace{1mm}
\State $\beta_{k+1} = \dfrac{1+\sqrt{1+4\beta_{k}^2}}{2}, \quad \gamma_k = \dfrac{1-\beta_k}{\beta_{k+1}}$ 
\State $y_{k+1} = x_k - \zeta_k G^{f_{\mu_{k+1}}}_{\zeta_k h}(x_k)$, \, $\zeta_k = \mu_{k+1}$ \vspace{1mm}
\State $x_{k+1} = (1-\gamma_k)y_{k+1} + \gamma_k y_k$ 
\end{algorithmic}
\end{algorithm}
Our first result quantifies the performance of Algorithm~\ref{alg:Algorithm1} after $T$~iterations for a specific choice of adaptive smoothing sequence $\big(\mu_k\big)_{k \ge 0}$.
\begin{lemma}[\small{Optimality gap in adaptive smoothing}]\label{prop:Alg1_bound}
Consider the optimization problem in \eqref{main problem} and Algorithm \ref{alg:Algorithm1} with adaptive smoothing variable {$ \mu_{k} = \max \Big\{\mu_{k-1}\Big(\dfrac{3\beta_{k}^2}{\beta_{k-1}^2}-1\Big)^{-1}, c \Big\}$}, for some $c \geq 0$. Then, after $T$ iterations, we have
\begin{align}\label{eq:prop:bound}
 & F(y_{T+1}) - F^\star \le \dfrac{L_{f}^2 \mu_{T+1}}{2} + \dfrac{E}{2 \beta_{T}^2\mu_{T+1}},
\end{align}
where the constant $E$ is 
\begin{align*}
    E = \|u_1\|^2 + \zeta_0\beta_{0}^2 \delta_{1} + \Big(1 - {\mu_1 \over \mu_0}\Big)\dfrac{3\mu_{0}}{2} L_{f}^2 \zeta_{0} \beta_{0}^2,
\end{align*} 
with $\delta_1 = f_{\mu_1}(y_1) + h(y_1) - F^\star$ and $u_1 = \beta_1 x_1 - (\beta_1 - 1)y_1 - x^\star$.
\end{lemma}
Lemma~\ref{prop:Alg1_bound} serves as a basis for different types of convergence results (global and local) in this study. {We note that, the right-hand side of the bound \eqref{eq:prop:bound} consists of two terms that are dependent on the adaptive smoothing parameter $\mu_k$ and can be used to control the optimality gap.} The proof of Lemma \ref{prop:Alg1_bound} builds on two preparatory lemmas, which we relegate to Section \ref{thechnicalproof} to improve the flow of the paper. 
Before proceeding with the main results concerning the convergence of Algorithm~\ref{alg:Algorithm1}, we first provide a lemma that sheds light on the behavior of the smoothing parameter~$\mu_k$ that we use in the algorithm.

\begin{lemma}[Smoothing parameter convergence]
\label{lem:rate}
    {The sequence generated by the first part of $\mu_k$ in Lemma \ref{prop:Alg1_bound}, i.e., $\mu_{k-1}\big({3\beta_{k}^2}/{\beta_{k-1}^2}-1\big)^{-1}$, exhibits a transient behavior at the order of~$\mathcal{O}(1/k^2)$ and an asymptotic linear rate of~$\mathcal{O}(e^{-k})$.}
\end{lemma}
The proof of Lemma~\ref{lem:rate} is rather straightforward and we defer it to Section \ref{thechnicalproof}. The convergence behavior of $\mu_k$ helps the global and local convergence guarantees of Algorithm~\ref{alg:Algorithm1}. In fact, our first result is to show that Algorithm~\ref{alg:Algorithm1} enjoys a global optimal convergence rate of $\mathcal{O}(1/\varepsilon)$ when the smoothing rule of~$\mu_k$ is uniformly lower bounded away from zero. 
\begin{theorem}[Global sublinear convergence]
\label{thm:global convegence}
Suppose the sequence of the smoothing parameters $\big(\mu_k\big)_{k\ge 0}$ is defined as in Lemma \ref{prop:Alg1_bound} where $c = {\varepsilon}/{L_f^2}$. Then, the outcome of Algorithm \ref{alg:Algorithm1} after $T$ iterations satisfies $F(y_T) - F^* \leq \varepsilon$ if $T \geq {2L_f\sqrt{E}}/{\varepsilon}$.
\end{theorem}
\begin{customproof}[Proof of Theorem~\ref{thm:global convegence}]
The proof builds on the assertion of Lemma~\ref{prop:Alg1_bound}. From Lemma~\ref{lem:rate}, we know that the first phase of $\mu_{k}$ is decreasing with the rate of at least~$\mathcal{O}(1/k^2)$. Therefore, the smoothing rule~$\mu_k$ also converges to ${\epsilon}/{L_{f}^2}$ with the similar rate of at least~$\mathcal{O}(1/k^2)$. It then suffices to determine the minimum number of iterations $T$ so that $\mu_{T+1}$ reaches ${\epsilon}/{L_f^2}$. To this end, note that since $\beta_k \ge k/2$ for all $k$, we have:
\begin{align*}
    F(y_{T+1}) - F^\star \le \dfrac{\varepsilon}{2} + \dfrac{2L_{f}^2 E}{{T}^2 \varepsilon}.
\end{align*} 
Hence, the $\epsilon$ precision is guaranteed if $T$ greater than ${2L_f\sqrt{E}}/{\varepsilon}$ which concludes the desired assertion.
\end{customproof}
We note that the global sublinear rate provided in Theorem~\ref{thm:global convegence} is worst-case optimal for the class of nonsmooth functions~\eqref{main problem}~\cite{nesterov2005smooth,nesterov2007smoothing}. Another key message of this study is that under a certain condition over the nonsmoothness, Algorithm~\ref{alg:Algorithm1} can achieve a linear, but local, convergence rate when the smoothing parameter follows the sequences in Lemma \ref{lem:rate} ($c = 0$ in Lemma \ref{prop:Alg1_bound}). {This two-phase behavior is evident in the smoothing rule as elucidated in Lemma \ref{lem:rate}.}

Let us remind that the stepsize in Algorithm~\ref{alg:Algorithm1} is dictated by the smoothing parameter (i.e., $\zeta_k = \mu_{k+1}$). It is surprising to have an exponentially fast diminishing stepsize in optimization algorithms. However, the following result addresses this phenomenon and demonstrates that, under a specific {$\infty$-locally strong convexity} condition, such a rate facilitates local linear convergence. 

\begin{theorem}[Local linear convergence]\label{thm:linear}
Consider the optimization problem in~\eqref{main problem} where the hypotheses of Lemma~\ref{prop:Alg1_bound} hold with $c = 0$. In addition, suppose the nonsmooth term~$f$ is {$\infty$-locally strongly convex} at the solution~$x^*$ of \eqref{main problem}, i.e., for any $\rho > 0$ the function~$f(x) - \rho \|x\|^2$ is locally convex in a neighborhood containing~$x^*$. {Then, for any initial condition $\mu_0 > 0$, the sequence~$\left(y_k\right)_{k \ge 0}$ generated by Algorithm~\ref{alg:Algorithm1} converges locally linearly to $x^*$, i.e., there exist $\varepsilon > 0$ and $\alpha \in (0,1)$ such that for all $x_0 \in \mathbb{B}_\varepsilon(x^*)$ there exists $k_0$ where for all $k \ge k_0$, we have $F(y_k) - F^* \le \alpha^{(k-k_0)}(F(y_{k_0}) - F^*)$.}
\end{theorem}
Before proceeding with the proof of this theorem, let us note that the key feature leading to local linear convergence is the fact that when $f$ is {$\infty$-locally strongly convex}, the condition number of its Hessian~$\partial^2f_{\mu}$ at its optimizer is uniformly bounded for all sufficiently small $\mu$. This allows the first-order accelerated algorithm to converge locally linearly and independently of the smoothness parameter~$\mu_k$. {It is also worth noting that the linear convergence rate is uniform for all the initial conditions $x_0 \in \mathbb{B}_\varepsilon(x^*)$. Let us formally explain this as follows.} 
\begin{customproof}[Proof of Theorem~\ref{thm:linear}]
Let us define $F_{\mu_k}(x) := f_{\mu_k}(x) + h(x)$. Note that using the basic property of the Moreau envelope in Lemma~\ref{lem:Moreau}\ref{lem-bound}, we have the inequality
\begin{align}\label{inequality for lin-conv}
    F(y_k) - F(x^*) \leq \dfrac{\mu_k}{2}L_{f}^2 + F_{\mu_k}(y_k) - F_{\mu_k}^* 
\end{align}
where $F_{\mu_k}^* = \min_{x \in \R^n} F_{\mu_k}(x) = F_{\mu_k}(x^*_{\mu_k})$.
Thanks to Lemma~\ref{lem:rate}, we know that the term $\mu_k$, and as such the first term of the upper bound~\eqref{inequality for lin-conv}, converges to zero exponentially fast. Therefore, it suffices to show that $F_{\mu_k}(y_k) - F_{\mu_k}^*$ converges to zero linearly with a rate independently from $\mu_k$. To this end, we recall that the Moreau envelope is locally $\mathcal{C}^2$-differentiable~\cite{poliquin1996generalized}. From the recent work \cite{grimmer2023landscape}, the maximum and minimum Hessian eigenvalues of $f_{\mu_k}$ at the optimal point $x_{\mu_k}^*$ can be bounded by
\begin{align*}
    \lambda_{\max}\left(\dfrac{\partial^2}{\partial x^2}f_{\mu_k}(x_{\mu_k}^*)\right) &\leq \mu_k ^{-1},\\
    \lambda_{\min}\left(\dfrac{\partial^2}{\partial x^2}f_{\mu_k}(x_{\mu_k}^*)\right) &\geq \left(\rho^{-1} + \mu_k\right)^{-1}.
\end{align*}
where $\rho$ is the strongly convex parameter of $f_{\mu_k}$ (i.e., $f_{\mu_k}(x) - \rho\|x\|^2/2$ is convex). This observation yields the bound on the condition number~$\kappa := \lambda_{\max} / \lambda_{\min}$ of the Hessian as $1 \leq \kappa \leq 1+(\rho \mu_k)^{-1}$. Recall that thanks to the {$\infty$-locally strong convexity} feature of $f$, we know that the parameter $\rho$ can be arbitrarily high as we are closer to the solution~$x^*$. In other words, there exists a sufficiently small $\varepsilon > 0$ where we can choose $\rho$ large enough for all $x \in \mathbb{B}_\varepsilon (x^*)$. Since the (accelerated)-proximal gradient descent methods converge linearly with the rate of $(1 - 1/\kappa)$~\cite{karimi2016linear}, we can deduce that for all sufficiently small enough $\mu_k$, Algorithm~\ref{alg:Algorithm1} enjoys linear convergence in a neighborhood of the solution~$x^*$.
\end{customproof}
Note that the $\infty$-local strong convexity of $f$ is essential for having locally linear convergence in Theorem 3.4. Namely, this property enables a linear decrease in the smoothing parameter while the strong convexity parameter, characterized by $\rho(x)$, grows sufficiently fast as we converge to the minimizer. Interestingly, this property holds for many popular nonsmooth terms commonly encountered in applications; see Section~\ref{simulation} for several such examples and also \cite{bello2022quadratic} (e.g., Theorem~3.1) for a detailed analysis of the explicit computation of such parameters. 


{We close this section with two remarks concerning the class of $\infty$-locally strongly convex functions and the initial value of the smoothing parameter in the proposed algorithm.}
\begin{remark}[\small{$\infty$-locally strong convexity vs. sharpness}]\label{re: infty local strongly convex}
    We note that the concept of ``$\infty$-locally strong convexity'' is closely related to the ``sharpness'' property \cite{Polyak1987}, but in a slightly weaker manner in the sense that the former is typically used locally, whereas the latter is posed globally. More specifically, a function $f$ is $\infty$-locally strongly convex with parameter $\rho^{\text{$\infty$-sc}}_\varepsilon$ (respectively, sharp with the parameter $\rho^{\text{sh}}$) when the following holds:
    \begin{align*} 
    &\text{$\infty$-locally strong convexity}: \,\frac{\rho^{\text{$\infty$-sc}}_\varepsilon}{2} \|x - x^*\|^2 \leq f(x) - f(x^*) \quad\,\, \forall \, x\in \mathbb{B} _{\varepsilon}(x^*), \, \text{where} \quad \rho^{\text{$\infty$-sc}}_\varepsilon \underset{\varepsilon \downarrow 0}{\longrightarrow} \infty, \\
    &\text{Sharpness}: \,\frac{\rho^{\text{sh}}}{2} \|x - x^*\| \leq f(x) - f(x^*) \quad \forall \, x.
    \end{align*}
\end{remark}

\begin{remark}[Avoid switching rule in Algorithm~\ref{alg:Algorithm1}]
\label{rem:initial}
In view of Theorem~\ref{thm:linear}, local linear convergence is achievable if we are sufficiently close to the global optimal point $x^*$. One can use Algorithm~\ref{alg:Algorithm1} with the lower-bounded smoothing rule (Theorem~\ref{thm:global convegence}) to converge to an arbitrarily close neighborhood of the optimal solution wherein linear convergence is guaranteed, and then continue with $c = 0$ in Algorithm~\ref{alg:Algorithm1} to converge to the desired solution with the faster rate anticipated in Theorem~\ref{thm:linear}. {To avoid this switching mechanism, a promising idea is to choose the initial value $\mu_0$ in a way that ensures a sufficient decrease in the objective function in the first phase behavior of $\mu_k$, and as such, guarantees reaching the desired neighborhood~$\mathbb{B}_\varepsilon (x^*)$ where the linear convergence is achievable. Note that the optimal value of $\mu_0$ typically depends on the initial error $\|x^* - x_0\|$, as we elaborate further in the next section.} 
\end{remark}
\subsection{Further discussion, limitation, and future direction}\label{further discussion section3}
In this part, we provide additional information and insights concerning the proposed smoothing technique, its limitations, and possible future directions.

\textbf{Further insights behind the adaptive smoothing rule:} The main motivation of the adaptive smoothness parameter is to exploit the possibility of having a larger smoothness parameter $\mu_k$ (and as such, optimizing smoother approximate function $f_\mu$), which leads to larger stepsizes and potentially faster convergence rate. However, any stepsize larger than $\epsilon/L_f^2$ can increase the Lyapunov function of Nesterov's accelerated algorithm, which is the key driving force behind our algorithm. This violation turns out to be dependent on the two subsequent algorithm momentum $\beta_k$ and $\beta_{k+1}$ and the two subsequent smoothness parameters $\mu_k$ and $\mu_{k+1}$.
This observation indeed leads to the smoothing rule in Lemma~\ref{prop:Alg1_bound}. Furthermore, as explained in Lemma~\ref{lem:rate}, the smoothing parameter $\mu_k$ has an initial decreasing rate of $\mathcal{O}(1/k^2)$, but asymptotically it converges to zero with an exponential rate. This behavior can explain the locally linear convergence of Algorithm~\ref{alg:Algorithm1} discussed in Theorem~\ref{thm:linear}.

\textbf{Extension to the smoothing rule:} It can be shown that the smoothing parameter $\mu_k$ in Lemma~\ref{prop:Alg1_bound} can be generalized to $\mu_{k} = \max \Big\{b\mu_{k-1}\big(\dfrac{b(a-1)+a}{a-1}\dfrac{\beta_{k}^2}{\beta_{k-1}^2}-1\big)^{-1}, c \Big\}$, where $a > 1$ and $b > 0$ are hyperparameters that control the behavior of the smoothing parameter $\mu_k$. These parameters can be selected to prevent Algorithm~\ref{alg:Algorithm1} from switching, thereby enabling locally linear convergence, as discussed in Remark~\ref{rem:initial}. The optimal tuning of these parameters remains unclear to the authors and is a promising direction for future research.

\textbf{Relation to prior works and the existing performance guarantees:} We note that our apriori theoretical results in Theorem~\ref{thm:global convegence} match the state of the art \cite{nesterov2005smooth} for the particular choice of algorithm parameters $(a,b) = (1,0)$, but do not improve the global performance. We wish also to note that this is not a rare precedent in optimization algorithm literature that a new algorithm numerically performs better than its formal apriori guarantees. For instance, the well-known FISTA algorithm proposed by~\cite{beck2009fast} has the same convergence bounds as in Nesterov's method \cite{nesterov2005smooth} when the latter is restricted to the proximal~setting.

\textbf{Limitation and future direction:} 
The primary limitation of the proposed method lies in the selection of the initial condition~$\mu_0$ (see Remark~\ref{rem:initial}). When $\mu_0$ is too small and $c = 0$, there is a risk that the algorithm fails to reach the region where linear convergence is attainable. This can result in stagnation of iterations due to the diminishing but still summable rate of decrease in the sequence of stepsizes~$\zeta_k = \mu_{k+1}$. Conversely, if $\mu_0$ is excessively large, the algorithm precision is compromised, as it takes longer to reach the linear convergence region. An optimal value for $\mu_0$ may depend on the initial error $\|x^* - x_0\|$ and the local characteristics of the functions involved in \eqref{main problem} at $x^*$~\cite{liang2017activity}. Investigating and analyzing these features are promising avenues for future research. 

Another future direction is concerned with the relation between the stepsize and smoothing parameters. In this work, the smoothing parameter dictates the stepsize. However, if we untangle this dependency, an adaptive stepsize may support the acceleration of the algorithms while an adaptive smoothing parameter enhances the precision. This adaptive stepsize-adaptive smoothing rule can be a promising research direction. 
\section{Technical proofs}\label{thechnicalproof}
In this section, we provide the theoretical proof for Section \ref{section3} and additional material supporting the technical as well as the numerical investigation of the paper. 
\subsection{Details of the Theoretical Analysis}\label{appendix}
\begin{figure}[!h]
     \centering           \includegraphics[width=0.35\textwidth]{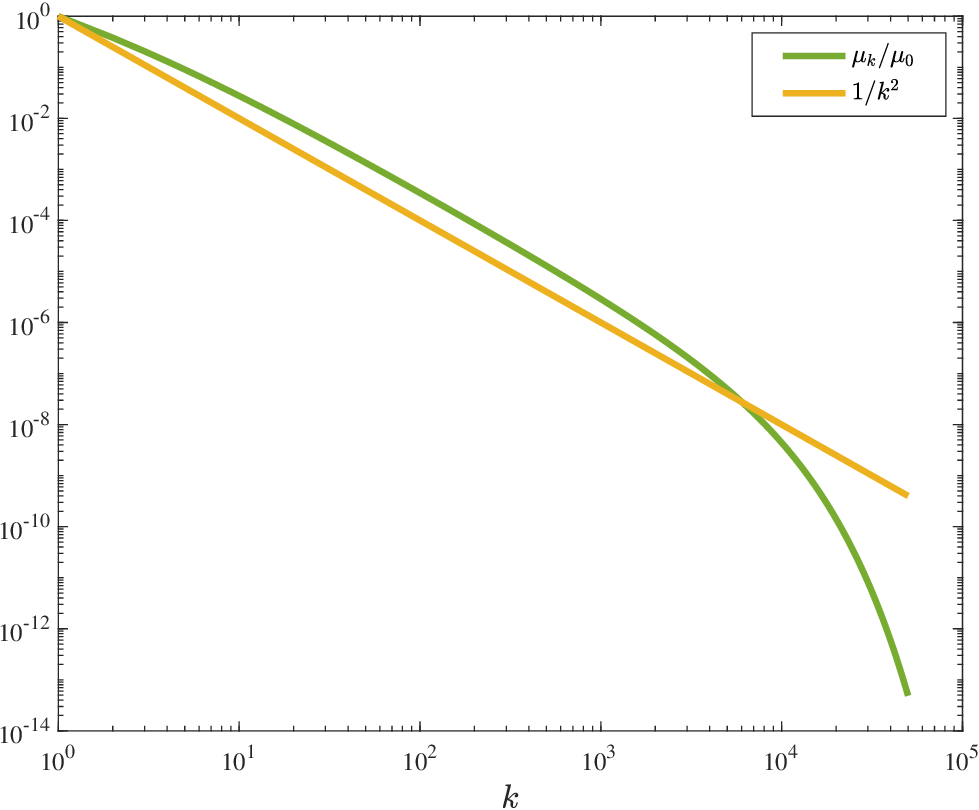}
           \caption{Convergence rate of $\mu_k$.}
           \label{smoothing parameter convergence rate}
\end{figure}

{Our proof for Lemma~\ref{prop:Alg1_bound} relies on a Lyapunov argument. To this end, we first proceed with two preliminary lemmas.}
\begin{lemma}[Gradient mapping]\label{Gradient mapping}
Let $G^{f}_{\zeta h}(x)$ be the gradient mapping defined in~\eqref{grad_mapping} where $\zeta$ is a positive scalar, $f(x)$ is smooth, and $h(x)$ is prox-friendly. Then, we have
\begin{align*}
    h\big(x-\zeta G^{f}_{\zeta h}(x)\big) \leq h(y) - \big\langle G^{f}_{\zeta h}(x) - \nabla f(x), y - \big(x - \zeta G^{f}_{\zeta h}(x)\big) \big\rangle, \quad \forall x,y \in \mathbb{R}^d.
\end{align*}
\end{lemma}
\begin{customproof}[Proof of Lemma~\ref{Gradient mapping}]
    Defining $u = \text{prox}_{\zeta h}(w)$, we can write
    \begin{align*}
        u = \text{prox}_{\zeta h}(w) \Leftrightarrow u = \arg \min_{u} h(u) + \dfrac{1}{2\zeta} ||u-w||^2 \Leftrightarrow  0 \in \partial{h(u)}+\dfrac{1}{\zeta}(u-w)\Leftrightarrow w-u \in \zeta \partial{h(u)}   
        \end{align*}
    By defining $u := x-\zeta G^{f}_{\zeta h}(x)$ and $w := x-\zeta \nabla f(x)$, we have
        \begin{align*}
         \underbrace{x-\zeta G^{f}_{\zeta h}(x)}_{u} &= x - \zeta \dfrac{1}{\zeta}(x - \text{prox}_{\zeta h}(x - \zeta \nabla f(x))) = \\
        &\text{prox}_{\zeta h}(\underbrace{x - \zeta \nabla f(x)}_{w}) \Rightarrow G^{f}_{\zeta h}(x) - \nabla f(x) \in \partial{h(x-\zeta G^{f}_{\zeta h}(x))}.
        \end{align*}
    Using the convexity of $h$, we can write
        \begin{align*}
        h(x-\zeta G^{f}_{\zeta h}(x)) \leq  h(y) - \langle G^{f}_{\zeta h}(x) - \nabla f(x), y - (x - \zeta G^{f}_{\zeta h}(x)) \rangle. 
        \end{align*}
\end{customproof}
Lemma~\ref{Gradient mapping} is an inherent property of the gradient mapping and essentially represents a convex inequality that is particularly helpful to control the increment of the original function~$F$ in \eqref{main problem}.
\begin{lemma}[Increment bound]\label{lem:increment}
Suppose function~$f$ is prox-friendly and $f_\mu$ is the smooth approximation~\eqref{smooth appr}. Considering the update $y_{k+1} = x_k - \zeta_k G^{f_{\mu_{k+1}}}_{\zeta_k h}(x_k)$, for any $z \in \R^d$ we have
\begin{align}\label{function reduction}
f_{\mu_{k+1}}(y_{k+1}) - f(z) + h(y_{k+1}) - h(z) \leq - \dfrac{1}{2\zeta_k}\|y_{k+1} - x_k\|^2 - \dfrac{1}{\zeta_k} \big\langle y_{k+1} - x_k , x_k - z \big\rangle
\end{align}
\end{lemma}
\begin{customproof}[Proof of Lemma~\ref{lem:increment}]
    By using the uniform boundedness of $f_{\mu}$ and the definition of convexity, we have
    \begin{align}
        f_{\mu_{k+1}}(y_{k+1}) &- f(z) + h(y_{k+1}) - h(z)\stackrel{\ref{(i)}}{\leq} f_{\mu_{k+1}}(x_k - \zeta_k G^{f_{\mu_{k+1}}}_{\zeta_k h}(x_k)) - f_{\mu_{k+1}}(z)  \nonumber\\
        &+ h(x_k - \zeta_k  G^{f_{\mu_{k+1}}}_{\zeta_k h}(x_k)) - h(z)\leq  f_{\mu_{k+1}}(x_k - \zeta_k G^{f_{\mu_{k+1}}}_{\zeta_k h}(x_k)) - f_{\mu_{k+1}}(x_{k}) \nonumber \\
        &+ \langle \nabla f_{\mu_{k+1}}(x_{k}), x_k - z \rangle -\underbrace{\langle G^{f_{\mu_{k+1}}}_{\zeta_k h}(x_k) - \nabla f_{\mu_{k+1}}(x_k), z - (x_k - \zeta_k G^{f_{\mu_{k+1}}}_{\zeta_k h}(x_k)) \rangle}_{\textbf{Lemma}\,\,\, \ref{Gradient mapping}} \nonumber\\
        &  \leq \langle \nabla  f_{\mu_{k+1}}(x_k), x_k - \zeta_k G^{f_{\mu_{k+1}}}_{\zeta_k h}(x_k) - x_k \rangle + \dfrac{\zeta_k^2}{2\mu_{k+1}}\|G^{f_{\mu_{k+1}}}_{\zeta_k h}(x_k)\|^2 + \langle G^{f_{\mu_{k+1}}}_{\zeta_k h}(x_k), y_{k+1} - z\rangle  \nonumber \\
        &  + \langle \nabla f_{\mu_{k+1}}(x_k), \zeta_k G^{f_{\mu_{k+1}}}_{\zeta_k h}(x_k) \rangle 
        \label{aux-eq3}
        \end{align}
    The last inequality is valid as a result of the smoothness property of the function $f_{\mu_{k+1}}$ \ref{lem-smooth}. By adding and subtracting $\dfrac{1}{\zeta_k} \langle y_{k+1} - x_k, x_k \rangle$ in \eqref{aux-eq3}, we obtain
        \begin{align*}
        f_{\mu_{k+1}}(y_{k+1}) - f(z) + h(y_{k+1}) - h(z) &\leq \underbrace {\langle y_{k+1} - x_k, \dfrac{1}{2\mu_{k+1}}\|y_{k+1} - x_k\|^2 -  \dfrac{1}{2\mu_{k+1}}\|y_{k+1} - x_k\|^2 \rangle}_{= 0}\\
        & - \dfrac{1}{2\zeta_k}\|y_{k+1} - x_k\|^2 - \dfrac{1}{\zeta_k}\langle y_{k+1} - x_k, x_k - z \rangle.
        \end{align*} 
\end{customproof}
Lemma~\ref{lem:increment} plays a key role in the proof of Lemma~\ref{prop:Alg1_bound}. The increment bound~\eqref{function reduction} allows for the inclusion of a momentum term that emerges in acceleration.
We are now in a position to prove Lemma~\ref{prop:Alg1_bound}. 
\begin{customproof}[Proof of Lemma~\ref{prop:Alg1_bound}]
{Here, we present the proof for the general selection of the smoothing parameter, as discussed in Section~\ref{further discussion section3}}. In the first step, by applying~\eqref{function reduction} and Lemma 2.6 in \cite{boct2015variable} for two cases of~$z = y_{k}$ and $z = x^*$ to arrive~at
\begin{subequations}
\begin{align}
     f_{\mu_{k+1}}(y_{k+1}) &- f_{\mu_{k}}(y_{k}) - \dfrac{(\mu_k - \mu_{k+1})}{2} L_{f}^2 + h(y_{k+1}) 
    - h(y_k) \stackrel{\text{\ref{(i)}}}{\leq} f_{\mu_{k+1}}(y_{k+1}) - f(y_k)\nonumber \\
    & + h(y_{k+1}) - h(y_k) 
    \leq - \dfrac{1}{2\zeta_k} \|y_{k+1} - x_k\|^2 - \dfrac{1}{\zeta_k}\langle y_{k+1} - x_k,  x_k - y_k \rangle. \label{eq4}\\
   f_{\mu_{k+1}}(y_{k+1}) &- f^* + h(y_{k+1}) - h^*  \leq 
   - \dfrac{1}{2\zeta_k} \|y_{k+1} - x_k\|^2 - \dfrac{1}{\zeta_k}\langle y_{k+1} - x_k, x_k - x^* \rangle. \label{eq5}
\end{align}  
\end{subequations}
Let us define~$\delta_k := f_{\mu_{k}}(y_{k}) + h(y_k) - f^* - h^*$. Then, multiplying \eqref{eq4} by $(\beta_k - 1)$ and adding the two sides of the inequality to \eqref{eq5} yields 
\begin{align*}
    &\beta_k \delta_{k+1} - (\beta_k - 1)\delta_k - \dfrac{(\mu_k - \mu_{k+1})}{2} L_{f}^2 (\beta_k - 1) \leq\\
    &- \dfrac{\beta_k}{2\zeta_k}\|y_{k+1}- x_k\|^2 -\dfrac{1}{\zeta_k}\langle y_{k+1} - x_k, \beta_k x_k - (\beta_k - 1)y_{k} - x^* \rangle. 
\end{align*}
Multiplying the above inequality by $\zeta_k \beta_k$ and considering $\beta_{k-1}^2 := \beta_{k}^2 - \beta_{k}$ and $\zeta_k \leq \zeta_{k-1}$, we~have
\begin{align}\label{eq6}
\zeta_k \beta_{k}^2 \delta_{k+1} - \zeta_{k-1} \beta_{k-1}^2 \delta_{k} &- \dfrac{(\mu_k - \mu_{k+1})}{2} L_{f}^2 \zeta_k \beta_{k-1}^2 \leq - \dfrac{1}{2} \Big(\|\beta_k (y_{k+1} - x_k)\|^2 \nonumber \\
&+ 2\beta_k \langle y_{k+1} - x_k, \beta_k x_k - (\beta_k - 1)y_{k} - x^*\rangle\Big)   
\end{align}
The right-hand side of~\eqref{eq6} can be equivalently written as
\begin{align}\label{eq7}
&\|\beta_k (y_{k+1} - x_k)\|^2 + 2\beta_k \langle y_{k+1} - x_k, \beta_k x_k - (\beta_k - 1)y_{k} - x^* \rangle = \nonumber \\
& \|\beta_k y_{k+1} - (\beta_k - 1)y_{k} - x^*\|^2 - \|\beta_k x_k - (\beta_k - 1)y_{k} - x^*\|^2. 
\end{align}
By substituting \eqref{eq7} into \eqref{eq6} and by rearranging the inequality, we have
\begin{align}\label{eq8}
\zeta_k \beta_{k}^2 \delta_{k+1} - \zeta_{k-1}\beta_{k-1}^2 \delta_{k} - \dfrac{(\mu_k - \mu_{k+1})}{2} L_{f}^2 \zeta_k \beta_{k-1}^2 &\leq - \dfrac{1}{2} \Big(\|\beta_k y_{k+1} - (\beta_k - 1)y_{k} - x^*\|^2 \nonumber \\
&- \|\beta_k x_k - (\beta_k - 1)y_{k} - x^*\|^2\Big)
\end{align}
Using the update rule of $x_{k+1}$ on the right-hand side of \eqref{eq8} reduces to 
\begin{align}\label{eq9}
\beta_k y_{k+1} - (\beta_k - 1)y_{k} - x^* = \beta_{k+1} x_{k+1} - (\beta_{k+1} - 1)y_{k+1} - x^*  
\end{align}
 which is equivalent to
\begin{equation*}
x_{k+1}=\dfrac{(-1 + \beta_k + \beta_{k+1})}{\beta_{k+1}}y_{k+1} + \dfrac{1 - \beta_k}{\beta_{k+1}}y_{k},   
\end{equation*}
By combining \eqref{eq8} and \eqref{eq9} with $u_k = \beta_k x_k - (\beta_k - 1) y_{k} - x^*$, we obtain
\begin{align}\label{eq10}
\zeta_k \beta_{k}^2 \delta_{k+1} - \zeta_{k-1} \beta_{k-1}^2 \delta_k  - \dfrac{(\mu_k - \mu_{k+1})}{2} L_{f}^2 \zeta_k \beta_{k-1}^2 \leq \dfrac{1}{2}\Big(\|u_k\|^2 - \|u_{k+1}\|^2\Big),
\end{align}
By defining $\omega_k := \mu_{k+1} / \mu_{k}$, we can rewrite \eqref{eq10} as
\begin{align}\label{eq10-2}
    \zeta_k \beta_{k}^2 \delta_{k+1} - \zeta_{k-1} \beta_{k-1}^2 \delta_k  - \dfrac{\mu_k(1 - \omega_k)}{2} L_{f}^2 \zeta_k \beta_{k-1}^2 \leq \dfrac{1}{2}\Big(\|u_k\|^2 - \|u_{k+1}\|^2\Big),
\end{align}
where $0 < \omega_k \leq 1$. Now, we consider two cases. First, we assume that $\mu_k$ strictly decreases in each iteration. To enforce $\mu_k$ being strictly decreasing we impose the monotonicity condition~$\mu_{k} < \mu_{k-1}$. Then, 
\begin{align}\label{eq11}
 - a \dfrac{\mu_{k-1}}{2} L_{f}^2 \zeta_{k-1} \beta_{k-1}^2 + (a - 1) \dfrac{\mu_k}{2} L_{f}^2 \zeta_k \beta_{k-1}^2 < - a \dfrac{\mu_k}{2} L_{f}^2 \zeta_k \beta_{k-1}^2 + (a - 1)\dfrac{\mu_k}{2} L_{f}^2 \zeta_k \beta_{k-1}^2 = - \dfrac{\mu_k}{2} L_{f}^2 \zeta_k \beta_{k-1}^2 
\end{align}
The inequality~\eqref{eq11} inspire us to define $\mu_k$ as 
\begin{align}\label{mu-def}
\mu_k = \left(\dfrac{b(a-1)+a}{a-1}\dfrac{\beta_k^2}{\beta_{k-1}^2} \mu_k - b \mu_{k-1} \right) \Rightarrow \, \mu_{k} = \dfrac{b\mu_{k-1}}{\dfrac{b(a-1)+a}{a-1}\dfrac{\beta_{k}^2}{\beta_{k-1}^2}-1}.
\end{align}
Using the definition of $\mu_k$ and its decreasing rate (Lemma \ref{lem:rate}), it can be easily shown that $1 - \omega_{k} \leq 1 - \omega_{k-1}$. Then, by substituting $\mu_k$ in the left-hand side of \eqref{eq11} and using \eqref{eq10-2}, we arrive at a Lyapunov-like inequality
\begin{align}\label{eq12}
\zeta_k \beta_{k}^2 \delta_{k+1} - \zeta_{k-1} \beta_{k-1}^2 \delta_k &-(1- \omega_{k-1})(ab-b+a) \dfrac{\mu_{k-1}}{2} L_{f}^2 \zeta_{k-1} \beta_{k-1}^2 \nonumber \\
&+(1 - \omega_{k})(ab-b+a) \dfrac{\mu_k}{2} L_{f}^2 \zeta_k \beta_{k}^2 \leq \dfrac{1}{2}\Big(\|u_k\|^2 - \|u_{k+1}\|^2\Big)
\end{align}
By summing up the inequalities in \eqref{eq12} from $ k = 1$ to $k = T$, one obtains
\begin{align} \label{Gap-appr-ineq}
\zeta_{T}\beta_{T}^2 \delta_{T+1} - \zeta_0\beta_{0}^2 \delta_{1} - (1 - \omega_0)(ab-b+a) \dfrac{\mu_{0}}{2} L_{f}^2 \zeta_{0} \beta_{0}^2 \leq \dfrac{1}{2}\|u_1\|^2 - \dfrac{1}{2}\|u_{T+1}\|^2 \leq \dfrac{1}{2}||u_1||^2,
\end{align}
which implies
\begin{align}\label{convergence eq}
&\delta_{T+1} \leq \dfrac{E}{2 \zeta_{T}\beta_{T}^2},\quad E = \|u_1\|^2 + \zeta_0\beta_{0}^2 \delta_{1} + (1 - \omega_0)(ab-b+a) \dfrac{\mu_{0}}{2} L_{f}^2 \zeta_{0} \beta_{0}^2.    
\end{align}
Second, we assume that the smoothing parameter $\mu_k$ strictly decreases by \eqref{mu-def} until it reaches to some a-priori value, and then we use the fixed $\mu_k$ afterward. By summing up the inequalities in \eqref{eq12} from $ k = 1$ to $k = K_\varepsilon$ (the iteration index that we fix $\mu_k$), one obtains
\begin{align} \label{2Gap-appr-ineq}
\zeta_{K_\varepsilon}\beta_{K_\varepsilon}^2 \delta_{K_\varepsilon+1} - \zeta_0\beta_{0}^2 \delta_{1} - (1 - \omega_0)(ab-b+a) \dfrac{\mu_{0}}{2} L_{f}^2 \zeta_{0} \beta_{0}^2 \leq \dfrac{1}{2}\|u_1\|^2 - \dfrac{1}{2}\|u_{K_\varepsilon+1}\|^2,
\end{align}
Now, by fixing $\mu_k$, we know that $\omega_k = 1$ for all $k \geq K_\varepsilon + 1$ and then by summing up the inequalities in \eqref{eq10-2} from $ k = K_\varepsilon + 1$ to $k = T$, one can obtain
\begin{align} \label{3Gap-appr-ineq}
\zeta_{T}\beta_{T}^2 \delta_{T+1} - \zeta_{K_\varepsilon}\beta_{K_\varepsilon}^2 \delta_{K_\varepsilon + 1} \leq \dfrac{1}{2}\|u_{K_\varepsilon + 1}\|^2 - \dfrac{1}{2}\|u_{T+1}\|^2,
\end{align}
Summing \eqref{2Gap-appr-ineq} and \eqref{3Gap-appr-ineq} yields the same inequality as \eqref{convergence eq}. Finally, by the definition of $\delta_{T+1}$, \ref{(i)}, and \eqref{convergence eq}, we then have
\begin{align}\label{final convergence eq}
F(y_{T+1}) - F^* \leq \dfrac{L_{f}^2}{2}\mu_{T+1} + \dfrac{E}{2 \zeta_{T}\beta_{T}^2}.    
\end{align}
{The result in Lemma~\ref{prop:Alg1_bound} can be easily derived by considering $a = 2$ and $b = 1$.}
\end{customproof}
\begin{customproof}[Proof of Lemma~\ref{lem:rate}]
{We present the proof for the general selection of the smoothing parameter, as discussed in Section~\ref{further discussion section3}}. Defining $\alpha_k := \big[\big({\beta_{k+1} \over {\beta_k}}\big)^2 + \big({a \over b(a-1)}\big)\big({\beta_{k+1} \over {\beta_k}}\big)^2 - {1 \over b} \big]^{-1}$, the first part of the smoothing parameter $\mu_k$ can be explicitly described as $\mu_{k+1} = \prod_{i \le k} \alpha_i$. To show $\mu_k \le {C_0 \over k^2}$ for some constant $C_0$, it suffices to show $\alpha_i \le \text{e}^{-2/i}$ for all $i \ge k_0$ where $k_0$ is a constant. This claim relies on the fact that this inequality implies the following:
\begin{equation*}
\mu_{k+1} = \prod_{i \le k} \alpha_i \le C_0\text{e}^{\sum_{i \le k} -2/i } \le C_0\text{e}^{-2\log(k)} \le {C_0 \over k^2},
\end{equation*}
where $C_0 = \prod_{i < k_0} \alpha_i$. To complete the proof, we show $\alpha_i \le \text{e}^{-2/i}$ for all sufficiently large $i$ by the following argument:\\
If $k_0 = 2\big(\log(1 + {1 \over b(a-1)})\big)^{-1} \le i$, then $\text{e}^{2/i} \le 1 + {1 \over b(a-1)}$, which consequently implies:
\begin{equation*}
    \text{e}^{2/i} \le \big({\beta_{i} \over {\beta_{i-1}}}\big)^2 + {1 \over b(a-1)} \Big( a\big({\beta_{i} \over {\beta_{i-1}}}\big)^2 - a + 1 \Big) = \alpha_i^{-1}.
\end{equation*}
In the last argument, we use the increasing property of $\beta_{i} \ge \beta_{i-1} > 0$ (i.e., $\big({\beta_{i} \over \beta_{i-1}}\big)^2 > 1$), which is a property of the momentum update (line 2 in Algorithm 1). Finally, to show the linear convergence, we show that $\alpha_k$ stays uniformly away from 1 as $k$ increases, ensuring a linear rate for large $k$. To this end, note that by the definition, we have
\begin{align*}
\alpha_k &:= \Big[\big({\beta_{k+1} \over {\beta_k}}\big)^2 + \Big({a \over b(a-1)}\Big)\big({\beta_{k+1} \over {\beta_k}}\big)^2 - {1 \over b} \Big]^{-1}  = \Big[\big({\beta_{k+1} \over {\beta_k}}\big)^2 \Big(1 + {a \over b(a-1)}\Big) - {1 \over b} \Big]^{-1} \\
& \le \Big[\Big(1 + {a \over b(a-1)}\Big) - {1 \over b} \Big]^{-1}  = \Big[1 + {1 \over b(a-1)} \Big]^{-1} = {b(a-1) \over b(a-1) + 1} < 1,
\end{align*}
where the second line inequality follows from the monotonicity of the momentum sequence $\beta_{k+1} \ge \beta_k$ and the fact that the coefficient $1 + a/b(a-1) > 0$ (the latter is also ensured by the feasible range of $a > 1$ and $b > 0$). This concludes that $\alpha_k$ is uniformly smaller than $1$ by the positive margin of $(b(a-1) + 1)^{-1}$.\\
For validation, Figure~\ref{smoothing parameter convergence rate} depicts the sequence in the first part of smoothing parameter $\mu_k$, as defined in Lemma \ref{lem:rate}. As illustrated in the figure, the smoothing parameter initially enjoys a convergence rate of $\mathcal{O}(1/k^2)$ but in the second phase, adopts an asymptotically linearly convergence rate. 
\end{customproof}
\section{Numerical experiments}\label{simulation}
We demonstrate the performance of Algorithm \ref{alg:Algorithm1} (with $c = 0$) on four popular classes of problems studied in the machine learning and control literature: (1) regression problems for both combination of the $\ell_1$ and $\ell_2$ norms borrowed from~\cite{yang2011alternating,thrampoulidis2015lasso}, (2) the \textit{MaxCut} problem belonging to the class of semidefinite programming from~\cite{sahin2019inexact}, (3) the Nuclear norm minimization problem and its application in model-free fault diagnosis from~\cite{avron2012efficient,noom2024proximal}, and (4) $\ell_1$-regularized model predictive control from \cite{annergren2012admm}. To evaluate our performance, we compare our proposed algorithm (Algorithm~\eqref{alg:Algorithm1}) with the following methods from the literature: (i) Sub-gradient descent (SGD)~\cite{nemirovskij1983problem}, the standard optimization algorithm, (ii)  Chambolle-Pock (CP)~\cite{chambolle2011first}, the state-of-the-art method for sparse regression, or the stochastic smoothing technique (stoch-smooth)~\cite{d2014stochastic}, a relatively recent method for semidefintie programming, (iii) Nesterov's smoothing (Nes-smooth)~\cite{nesterov2005smooth}, and (iv) Tran-Dinh (TD)~\cite{tran2017adaptive}, that is the closest in spirit to our proposed method. 

Note that {in light of Remark \ref{re: infty local strongly convex}, we use sharp convex functions in this study, and the simulations include examples that satisfy the sharpness property. An example of an $\infty$-locally strongly convex function is $\|\cdot\|_1 $. Using the definition of $\infty$-locally strongly convex functions and the sharpness properties of $\|\cdot\|_1$, there exists a neighborhood of $x^*$ where, for any $\rho > 0$, the function is $\rho$-strongly convex \cite{Polyak1987,liang2017activity,bello2022quadratic}. Several important prox-friendly functions fall into the category as well, including: $\|x\|_\infty$ and the Nuclear norm of matrices (see the tables "Prox Calculus Rules" and "Prox Computations" on pages 448 and 449 in \cite{beck2017first}). Additional examples, formed as combinations of such functions, are also discussed in the literature, for instance, \cite{liang2014local} considers the total variation, $\ell_\infty$-norm, $\ell_1-\ell_2$-norm, and the Nuclear norm. Indeed, the widespread applicability of these functions across various domains motivated us to include several numerical examples in this section. We would also like to highlight that the class of problems studied in this work encompasses several key applications in image processing, image reconstruction, system identification, and control. Specifically, in \cite{fu2006efficient}, the authors employ $\ell_2-\ell_1$ and $\ell_1-\ell_1$ optimization for image restoration. In \cite{chambolle2011first}, the TV-$\ell_1$ model is used for image denoising, while in \cite{wang2013dictionary}, the $\ell_1-\ell_1$ problem is applied to dictionary learning. Lastly, the $\ell_1$-Nuclear norm problem is employed in data-driven fault diagnosis control in \cite{noom2024proximal}.}
\subsection{Regression problems}\label{sec:regression}
We consider the $\ell_i$-$\ell_1$-regularized LASSO problem
\begin{equation}\label{regression problem}
F^* := \min\limits_{x \in \mathbb{R}^n}   \mathcal{R}(x) + \eta \|x\|_1,
\end{equation}
where  $\mathcal{R}(x) = \|Bx - b\|_i, \,\, i = \{1,2\}$, is a regularized function. the parameters~$B$ and $b$ are given and $\eta \geq 0$ is a weighting coefficient. Functions $f$ and $h$ in  \eqref{main problem} can be written as $f(x) := \mathcal{R}(x)$ and $h(x) := \eta \|x\|_1$, respectively. For the function~$f$, we use the smoothing approximation~\eqref{smooth appr}.

The test data is generated as follows: Matrix $B \in \mathbb{R}^{n\times n}$ is generated randomly using the standard Gaussian distribution $\mathcal{N}(0, 1)$. For simulation, we consider two cases as in~\cite{tran2017adaptive}. In the first case, we use uncorrelated data, while in the second case (reported in the Appendix~\ref{Supplemental Numerics}), we generate $B$ with $50\%$ correlated columns as $B(:,j+1) = 0.5B(:,j) + \textrm{randn}(:)$. The observed measurement vector $b$ is generated as $b := Bx^{\natural} + \mathcal{N}(0, 0.05)$, where $x^{\natural}$ is generated randomly using $\mathcal{N}(0, 1)$.  For the simulation, we set $n=100$ and the coefficient parameter $\eta$ is chosen as suggested in \cite{tran2017adaptive,bickel2009simultaneous}. All algorithms are initialized at the same point chosen randomly. In Nes-smooth method the stepsize is ${2\cdot10^{-3}}/{L_{f}^2\|B\|^2}$ for achieving desired error $\epsilon = 10^{-3}$ and the stepsizes in primal and dual updates in CP are ${1}/{\|B\|}$. In TD and the proposed method, the stepsize is $\zeta = \mu_k/\|B\|^2$. The convergence speed is significantly affected by the initial condition of $\mu$ in TD, and following the suggestion by \cite{tran2017adaptive}, we set $\mu_0 = \mu^* = {\|B\| \cdot \|x_0 - x^* \|}/{\sqrt{3L_{f}^2}}$ for the TD algorithm.

\textbullet \, \textbf{$\ell_1$-$\ell_1$-regularized LASSO.} Here, $i = 1$ for $\mathcal{R}(x)$ in \eqref{regression problem}. Then, $f$ can be written as $f(x) := \|Bx - b\|_1 = \max\limits_y \bigl\{ \langle B^\top y,x \rangle - \langle b,y \rangle: \|y\|_{\infty} \leq 1 \bigl\}$. Hence, we can smooth $f$ using the quadratic prox-function to obtain $f_{\mu}(x) := \max_{y}\Big\{ \langle B^\top y,x \rangle - \langle b,y \rangle - \dfrac{\mu}{2}\|y\|^2 : y \in \mathbb{B}_{\infty} \Big\}$. Using the analytical solution, one can observe that $\nabla f_{\mu}(x)= \mathrm{proj}_{\mathbb{B}_{\infty}}\left(\mu^{-1}(Bx - b)\right)$ and $L_{f}^2 = n$.

\textbullet \, \textbf{$\ell_2$-$\ell_1$-regularized LASSO.} We consider $i = 2$ for $\mathcal{R}(x)$ in \eqref{regression problem}. Functions $f$ can be written as $f(x) := \|Bx - b\|_2\ = \max\limits_y \bigl\{ \langle B^\top y,x \rangle - \langle b,y \rangle: \|y\|_{2} \leq 1 \bigl\}$. Hence, we similarly have $f_{\mu}(x) := \max_{y}\Big\{ \langle B^\top y,x \rangle - \langle b,y \rangle - \dfrac{\mu}{2}\|y\|^2 : y \in \mathbb{B}_{2} \Big\}$. In this case, we can  show that $\nabla f_{\mu}(x) =  \mathrm{proj}_{\mathbb{B}_{2}}\left(\mu^{-1}(Bx - b)\right)$, and $L_{f}^2 = 1$.

Figures \ref{fig1:sfig1} and \ref{fig1:sfig2} report the results for \eqref{regression problem} with $i = 1,2$ where all algorithms are initialized at the same point chosen randomly. As shown in Figures~\ref{fig1:sfig1} and \ref{fig1:sfig2}, the proposed algorithm outperforms other methods.
\subsection{\textit{MaxCut} (semidefinite programming)} \label{sec:maxcut}
An important application of the first-order smoothing technique is to solve semidefinite programming. This class can essentially be reduced to minimizing the maximum eigenvalue of a matrix as follows: 
\begin{equation}\label{max-eig}
\min\limits_{X \in \mathbb{S}^n}  \lambda_{\textbf{max}}(X)
\end{equation}
It is worth noting that the primal form of all semidefinite programs with a fixed trace can be rewritten as \eqref{max-eig}. As shown in \cite{nesterov2007smoothing}, the smooth approximation of $\lambda_{\textbf{max}}(X)$ can be written as $f_{\mu}(X) = \mu \log \left(\sum_{i = 1}^n \exp(\lambda_i(X)/\mu)\right)$, which is convex and twice differentiable with a gradient 
\begin{align*}
\nabla f_{\mu}(X) &=\\
&\left(\sum_{i = 1}^n \exp(\lambda_i(X)/\mu)\right)^{-1} {\sum_{i = 1}^n \exp(\lambda_i(X)/\mu) q_i q_i ^\top},
\end{align*}
where $q_i$ is the $i^{th}$ column of the unitary matrix $Q$ in the eigen-decomposition $Q \Sigma Q^T$ of $X$ and $\lambda_1 \geq . . . \geq \lambda_n$ are the eigenvalues of $X$. In addition, $f_{\mu}(X)$ fulfills the uniform boundedness inequality~\ref{(i)} with $L_{f}^2 = 2 \log n$ and the smoothness parameter~$ 1/\mu$ (cf. \ref{lem-smooth} in Lemma~\ref{lem:Moreau}).
The numerical performance of the smoothing technique is evaluated for the \textit{MaxCut} relaxation which is written in the primal form of the semidefinite program as follows:
\begin{align}\label{maxcut primal}
    \max\limits_{X \in \mathbb{S}^n} \textbf{Tr}(CX) \qquad \text{s.t.} \quad
    \text{diag}(X) = 1, \quad X \succcurlyeq 0 
\end{align}
Matrix $C$ is generated using the Wishart distribution with $C = {G^\top G}/{\|G\|_{2}^2}$, where $G$ is a standard Gaussian matrix. Here, we consider a regularized dual form of \eqref{maxcut primal} as suggested in \cite{helmberg2000spectral}
\begin{equation}\label{maxcut dual with penalty}
\min\limits_{y \in \mathbb{R}^n}  \lambda_{\textbf{max}}(C+\text{diag}(y)) - \langle \mathbf{1}, y \rangle + \eta \mathcal{R}(y),
\end{equation}
where $\eta$ is a regularization parameter. Figure~\ref{fig2:sfig1} and \ref{fig2:sfig2} show the simulation results for minimizing problem~\eqref{maxcut dual with penalty} with $n = 100$, $2000$ iterations, and two different choices of $\mathcal{R}(y)$. In the Nes-smooth method, the stepsize is $10^{-3}$ for achieving a desired error~$\epsilon = 10^{-3}$. The parameters required in the stoch-smooth method are chosen according to~\cite{d2014stochastic}. The initial condition for $\mu$ in the TD method is $\mu_0 = \mu^* = {\|x_0 - x^* \|}/{\sqrt{6 \log n}}$. Additional simulations for different $\eta$ are reported in the Appendix~\ref{Supplemental Numerics}.
\subsection{Nuclear norm minimization (Model-free Fault Diagnosis)} \label{sec:Nuclearnorm}
\begin{figure*}[h!] 
\centering
\begin{minipage}{0.32\textwidth}
  \includegraphics[width=\linewidth]{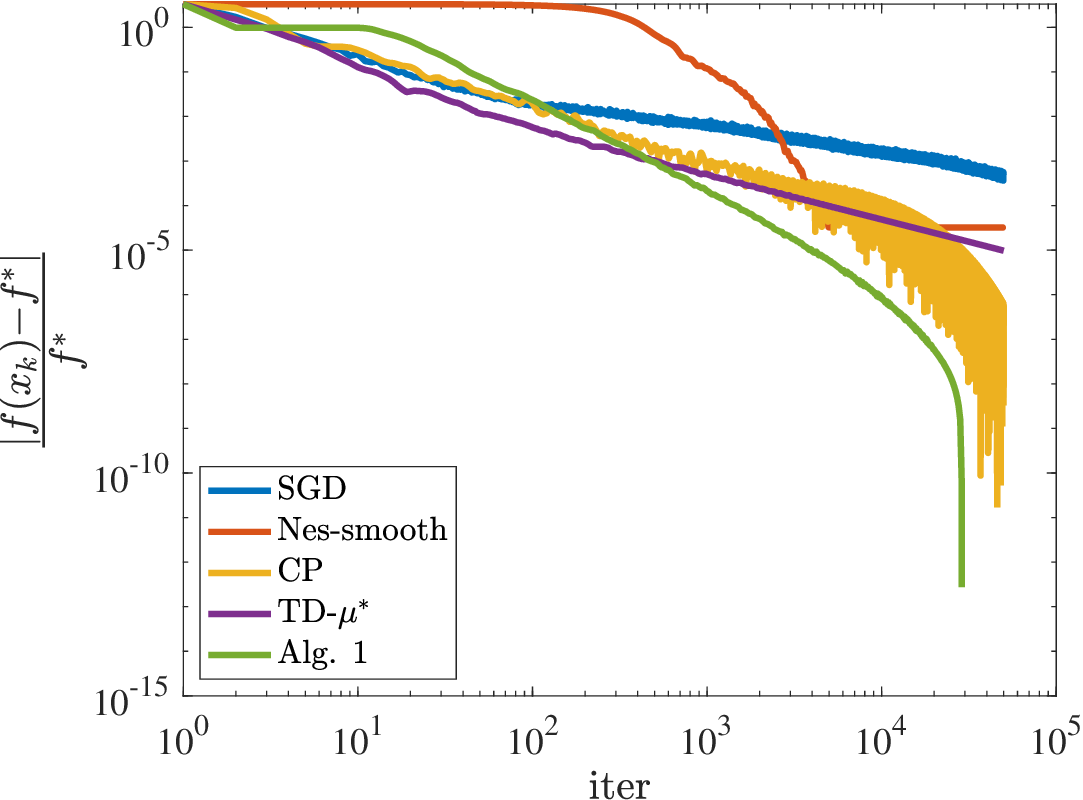}
  \subcaption{$\ell_1$-$\ell_1$-regularized LASSO~\eqref{regression problem}.}
  \label{fig1:sfig1}
\end{minipage}%
\hspace{-1mm}
\begin{minipage}{0.32\textwidth}
\includegraphics[width=\linewidth]{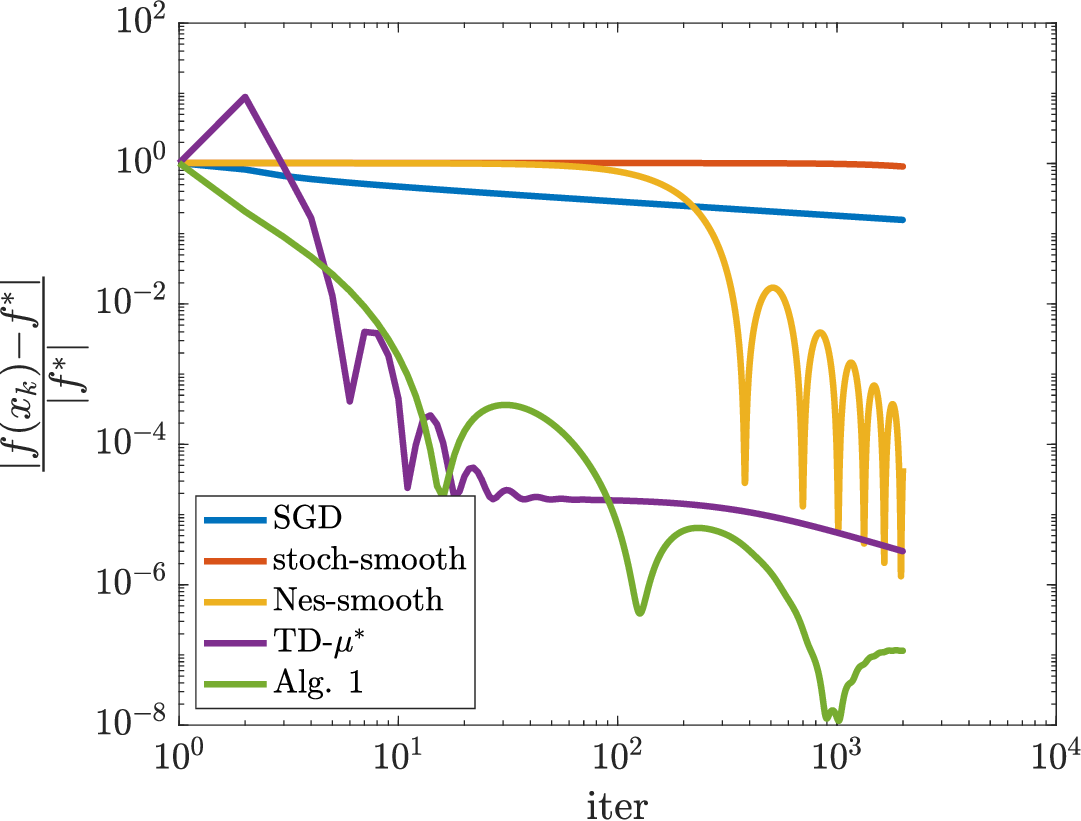}
        \subcaption{$\mathcal{R}(y) = \|y\|^2_{2}$ and $\eta = 0.05$~\eqref{maxcut dual with penalty}.}
  \label{fig2:sfig1}
\end{minipage}%
\hspace{-1mm}
\begin{minipage}{0.35\textwidth}
\includegraphics[width=\linewidth]{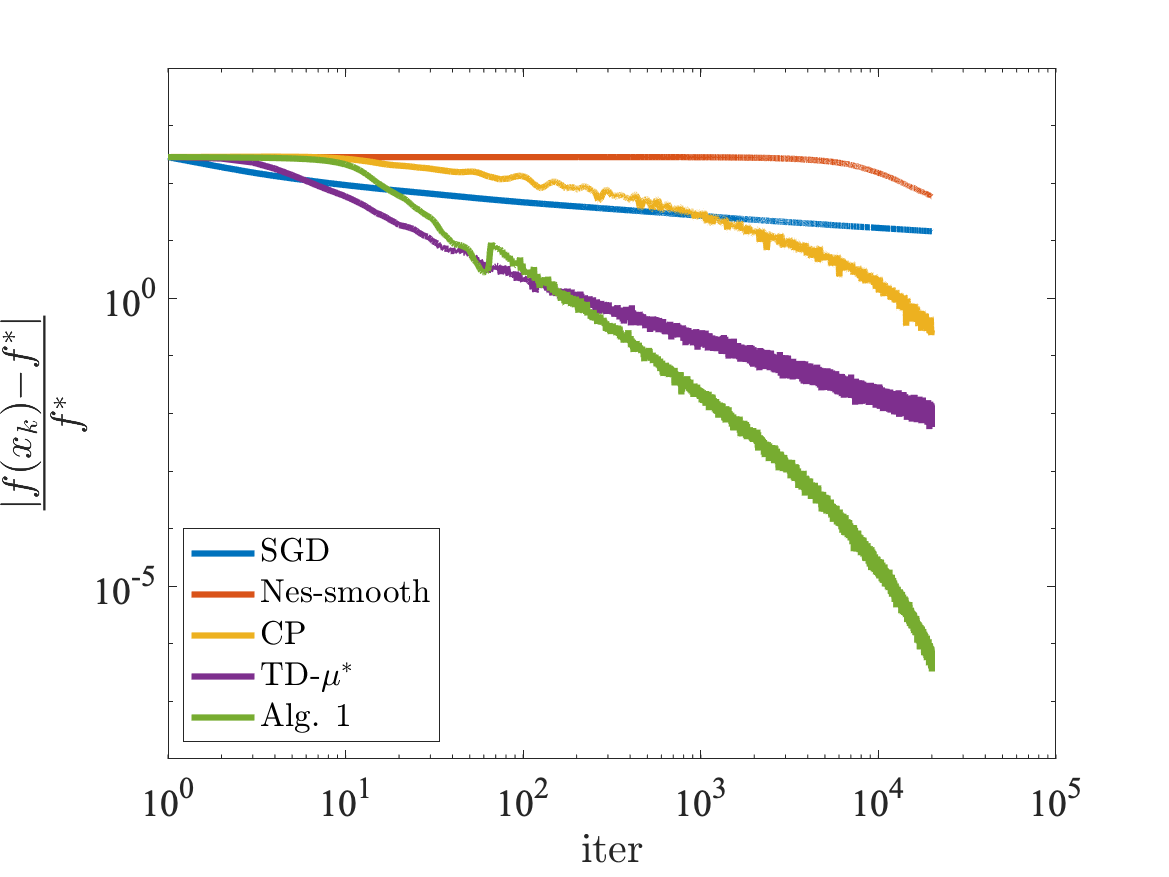}
       \subcaption{Fault diagnosis~\eqref{Nuclear-norm}.}
  \label{fig7:sfig1}
\end{minipage}
\begin{minipage}{0.32\textwidth}
  \includegraphics[width=\linewidth]{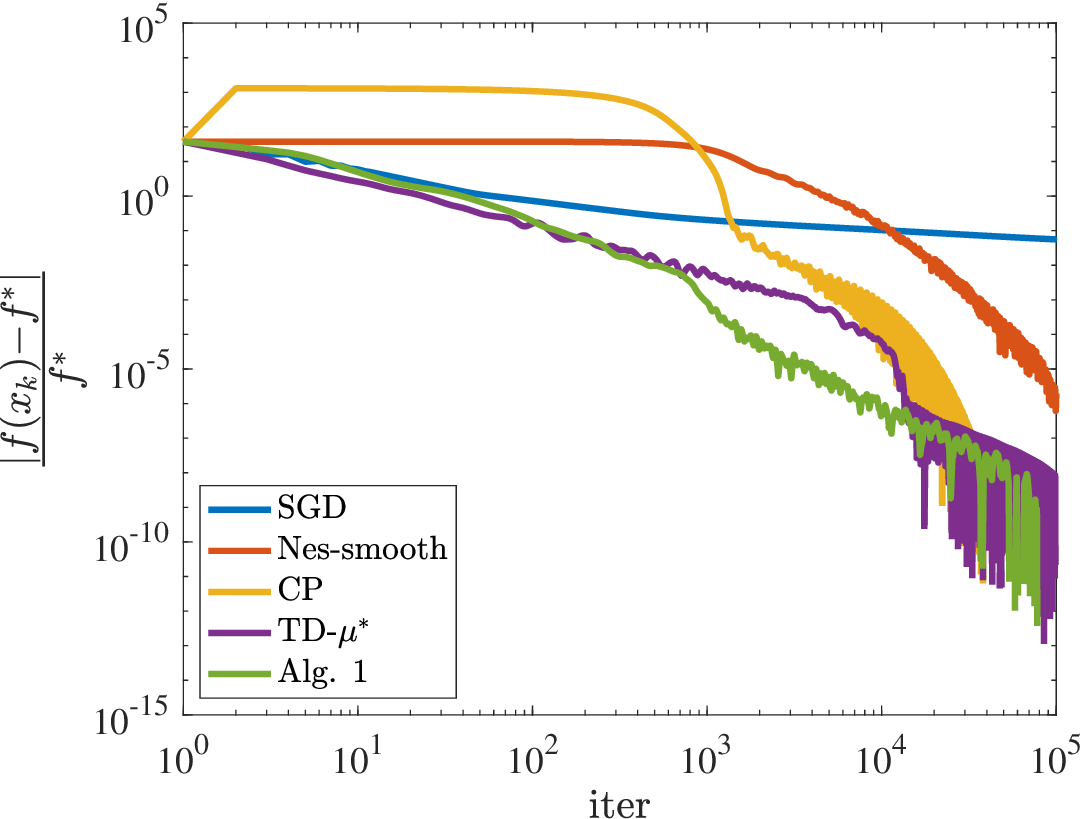}
      \subcaption{$\ell_2$-$\ell_1$-regularized LASSO~\eqref{regression problem}.}
  \label{fig1:sfig2}
\end{minipage}%
\hspace{-.10cm}
\begin{minipage}{0.32\textwidth}
\includegraphics[width=\linewidth]{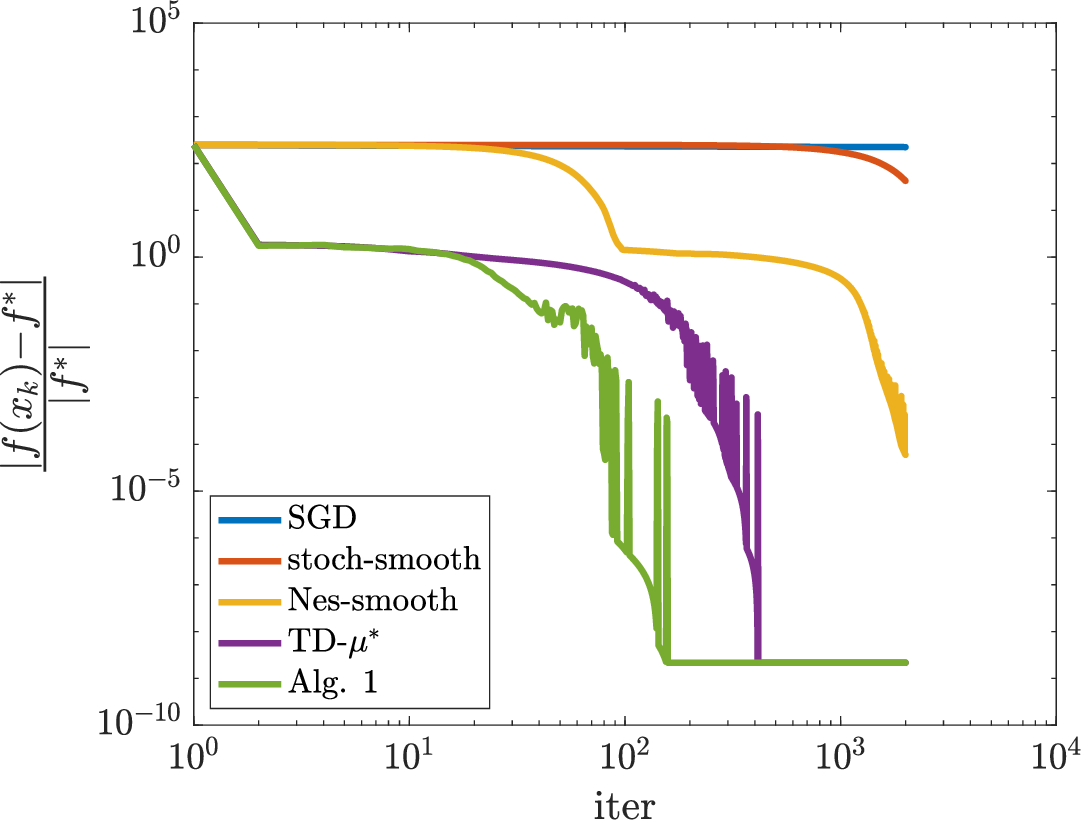}
      \subcaption{$\mathcal{R}(y) = \|y\|_1$ and $\eta = 1$~\eqref{maxcut dual with penalty}.}
  \label{fig2:sfig2}
\end{minipage}%
\hspace{-.10cm}
\begin{minipage}{0.32\textwidth}
\includegraphics[width=\linewidth]{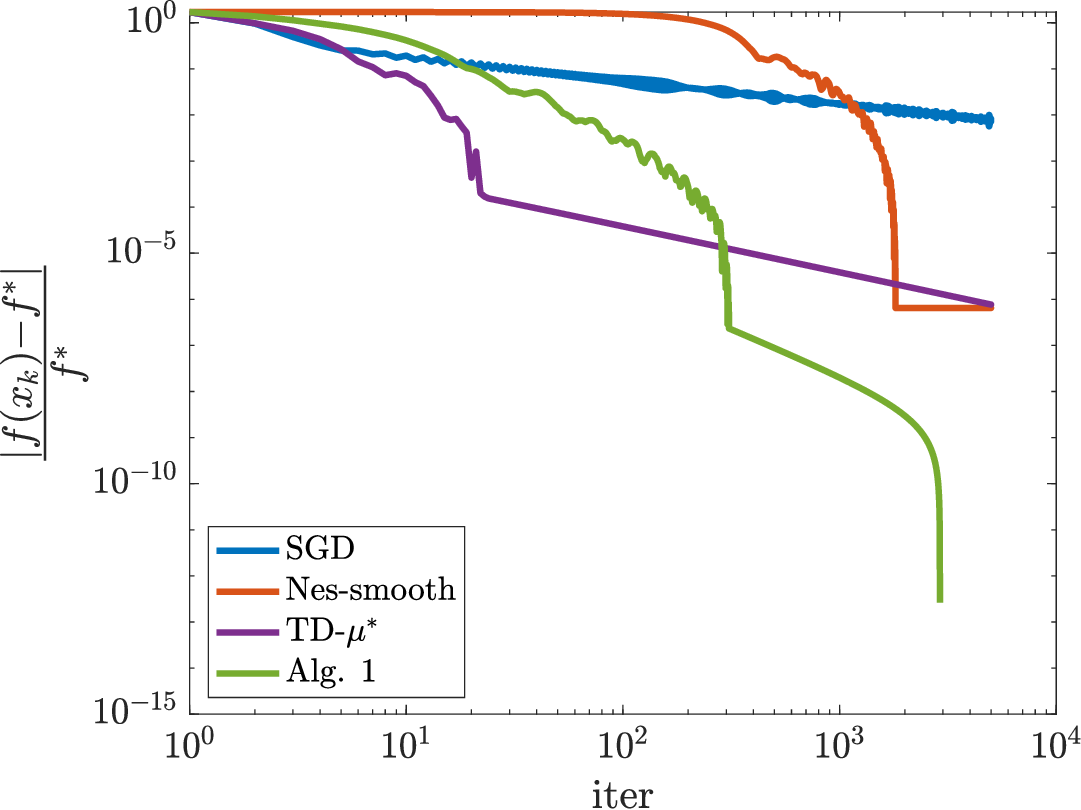}
    \subcaption{$\ell_1$-Regularized MPC~\eqref{regularized-mpc}.}
    \label{fig7:sfig2}
\end{minipage}
\caption{Numerical results for four problem classes. The first column shows regression problems with uncorrelated data (Subsection~\ref{sec:regression}); the second, the \textit{MaxCut} problem with different choices of $\mathcal{R}(y)$ (Subsection~\ref{sec:maxcut}); and the third, model-free fault diagnosis with Nuclear norm minimization (Subsection~\ref{sec:Nuclearnorm}) and $\ell_1$-Regularized MPC (Subsection~\ref{sec:mpc}).}
\label{fig:2x3_fullwidth}
\end{figure*}

The Nuclear norm regularization is an important problem emerging in several applications including matrix completion~\cite{avron2012efficient,jaggi2010simple}, compressed sensing~\cite{srebro2004maximum}, principal component analysis~\cite{guo2013robust}, system identification and fault diagnosis~\cite{sun2022system}. In this work, we consider model-free fault diagnosis as in \cite{noom2024proximal}, where the authors proposed a model-free, data-driven fault diagnosis approach that aims to identify the system and diagnose faults simultaneously, thereby eliminating the need for an extensive identification phase prior to fault detection. The proposed method reformulates the problem as a convex optimization problem involving the summation of nonsmooth terms. To enhance computational efficiency, proximal and splitting-type algorithms (similar to Nes-smooth, CP, and TD methods) are employed for online implementation. The problem formulation is expressed as follows (see equation (5) in \cite{noom2024proximal}):
\begin{align}\label{Nuclear-norm}
  \blank{-.3cm}  F^* := \min_{x, X_L, X_S}   \frac{1}{2}\|y - Hx\|^2 + \tau \|X_L\|_* + \lambda \|X_S\|_1.
\end{align}
Here, we consider the same buck converter example in \cite{noom2024proximal} with same parameters and solve the problem using different methods. Here, $X_L\in \mathbb{R}^{m_1 \times m_2}$ and $X_S\in \mathbb{R}^{n_1 \times n_2}$ are linear maps of $x$ that select and rearrange its elements into a low-rank matrix and a sparse vector, respectively. $H$ is a Toeplitz matrix corresponding to the input, output, and fault signals, $y$ denotes the output data, and $\tau \geq 0$ and $\lambda \geq 0$ are regularization parameters. The Nuclear norm of a matrix $X$, denoted $\|X\|_*$, is the sum of its singular values, equivalently the $\ell_1$-norm of the singular values. Both nonsmooth terms in \eqref{Nuclear-norm} are prox-friendly.  Therefore,  functions $f$ and $h$ in \eqref{main problem} can be expressed as $f(x) = \frac{1}{2}\|y - Hx\|^2 + \lambda \|X_S\|_1$ and $h(x) = \tau \|X_L\|_*$ and can be smoothed using \eqref{smooth appr}. \\
Figures~\ref{fig7:sfig1} depict the results of optimizing \eqref{Nuclear-norm}. All algorithms are initialized at the same randomly chosen point and run for 20,000 iterations. For the Nes-smooth method, the stepsize is set to ${2 \cdot 10^{-3}}/(\|H\|^2 + \lambda)$ to achieve a precision level of $\epsilon = 10^{-3}$. In the primal and dual updates of CP, the stepsizes are set to $1/(\|H\|^2 + \lambda)$, while in TD and the proposed method, the stepsize is determined based on $\mu_k / (\|H\|^2 + \lambda)$.
\subsection{$\ell_1$-regularized model predictive control} \label{sec:mpc}
The authors in \cite{annergren2012admm} consider the $\ell_1$-regularized MPC problem with horizon $H$ to reduce actuator activity in the quadruple tank system. They use a linearized state-space model and solve the following MPC problem \cite[Eq.~(1)]{johansson1999teaching}:
\begin{align}\label{regularized-mpc}
\min_{x,u} \quad  J(x,u) &:= \sum_{i=0}^{H-1} \Big( \| x_i - x_{\mathrm{ref},i} \|_Q^2 + \lambda \| \Delta u_i \|_1 \Big) + \| x_H - x_{\mathrm{ref},H} \|_Q^2 \\
\text{s.t.} \quad  x_{i+1} &= \textbf{A} x_i + \textbf{B} u_i, \quad \forall i \in \{0, \dots, H-1\}, \nonumber\\
\quad u_{\min} &\leq u_i \leq u_{\max}, \qquad \forall i \in \{0, \dots, H-1\},\nonumber \\
\Delta u_i &= u_i - u_{i-1}, \quad \forall i \in \{0, \dots, H-1\}.\nonumber
\end{align}
We use the same parameters as in \cite{annergren2012admm} for the system dynamics and solve problem \eqref{regularized-mpc} using different methods. The objective function can be represented by $f = J(x,u)$ and $h = {i}(\chi)$ in \eqref{main problem}, where $\chi$ is the feasible set of problem \eqref{regularized-mpc}, and ${i}(\chi)$ is its indicator function. Therefore, the nonsmooth $\ell_1$ term in $f$, which is prox-friendly, can be smoothed. Figure~\ref{fig7:sfig2} shows the results; as can be seen and mentioned in Remark \ref{rem:initial}, the proposed method exhibits two convergence phases, including an initial $\mathcal{O}(1/k^2)$ phase followed by linear convergence. Note that CP method cannot be applied to this problem because we have more than one nonsmooth terms in the objective function.

We provide additional numerical experiments using the real-world dataset LIBSVM \cite{chang2011libsvm}, as well as extensive simulations with increasing problem dimensions. These results are presented in the Appendix. Furthermore, we note that the ADMM method is suitable for some of the aforementioned problems, and additional simulations for comparison are provided in the Appendix as well. However, ADMM may fail to converge, converge very slowly, or be difficult to implement when multiple nonsmooth terms are present, particularly when they appear in a summation that cannot be easily split \cite{chen2016direct}.
\appendix
\section{Additional Numerical Results}\label{Supplemental Numerics}
This section presents additional numerical experiments showcasing the performance of the proposed method in the case of regression problem \eqref{regression problem} and \textit{MaxCut} \eqref{maxcut dual with penalty} and also the comparison with the ADMM and Algorithm \ref{alg:Algorithm1} with constant term $c = {\epsilon}/{L_{f}^2}$ (Theorem~\ref{thm:global convegence}, Alg. 1-2 in Figures).

\begin{figure}[h!]
\centering
\subfloat[$\ell_1$-$\ell_1$-regularized LASSO.]{\label{fig4:sfig1}\includegraphics[scale=0.32]{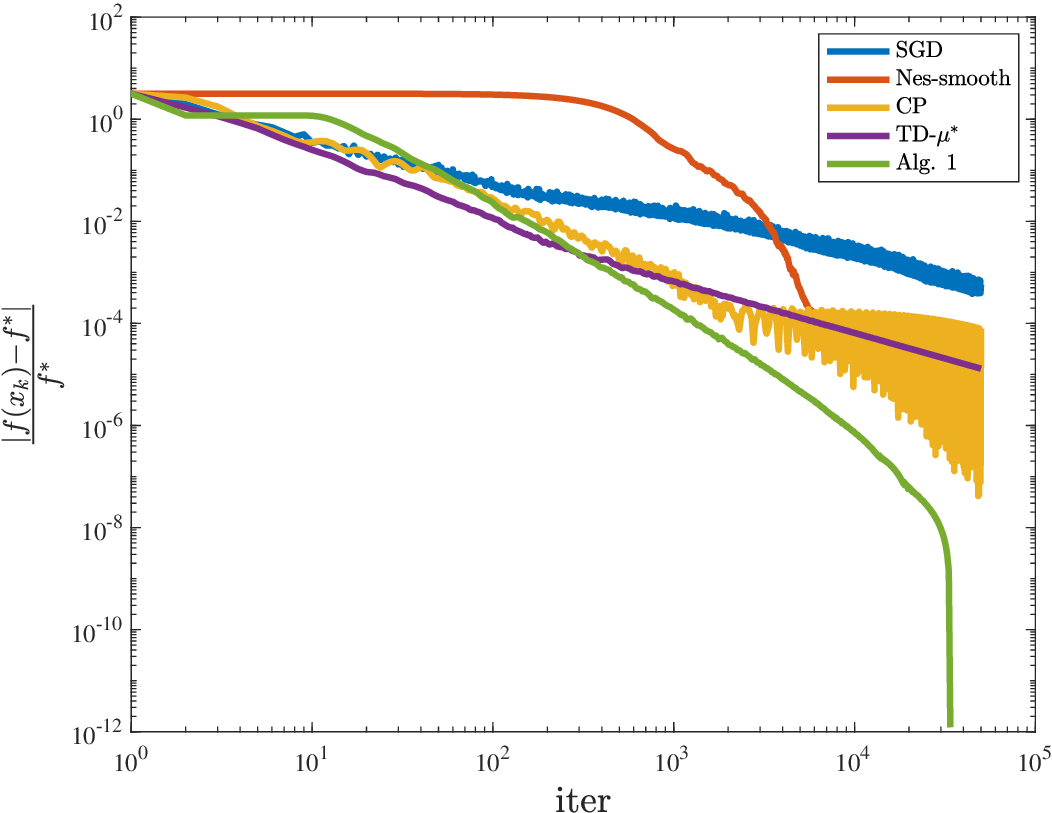}}
\hspace{8mm}
\subfloat[$\ell_2$-$\ell_1$-regularized LASSO.]{\label{fig4:sfig2}\includegraphics[scale=0.32]{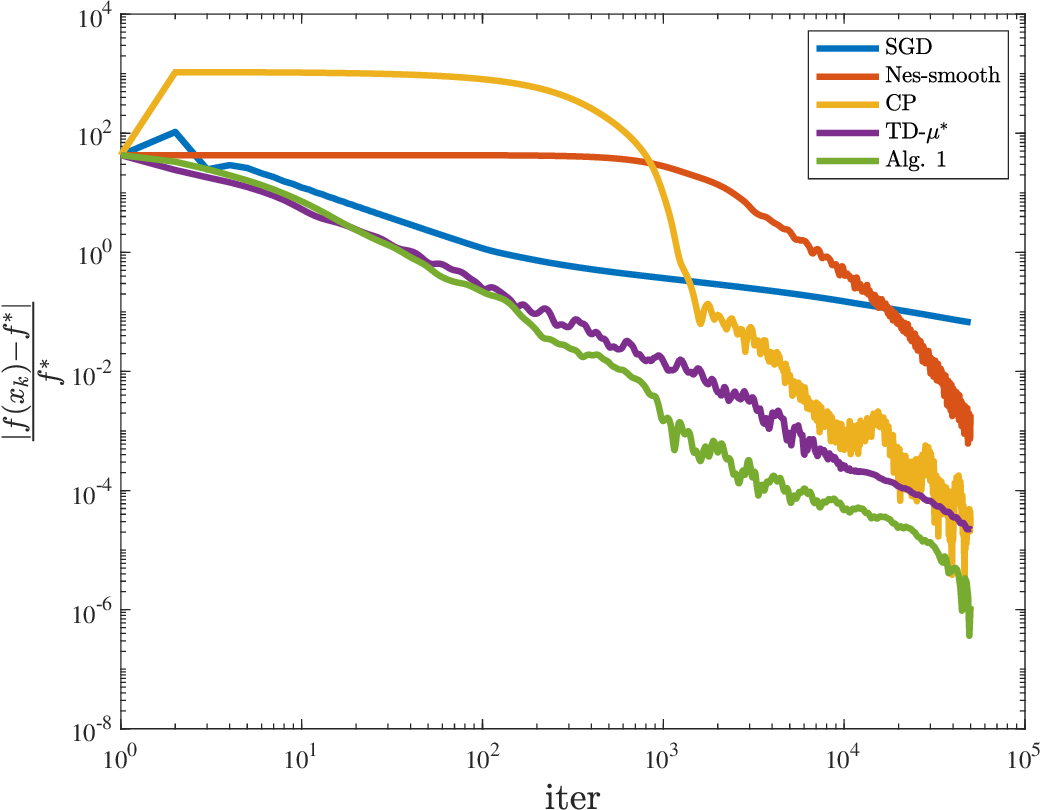}}\\
\caption{Numerical results for regression problems in Subsection~\ref{sec:regression} (correlated data).}
\label{LASSO-correlated}
\end{figure}

\textbullet \,\ Figure \ref{LASSO-correlated} contains plots for \eqref{regression problem} with $i = 1,2$ and correlated data where all algorithms are initialized at the same point chosen randomly.

\begin{figure}[h!]
\centering
\subfloat[$\mathcal{R}(y) = \|y\|^2_{2}$ and $\eta = 0.005$.]{\label{fig5:sfig1}\includegraphics[scale=0.32]{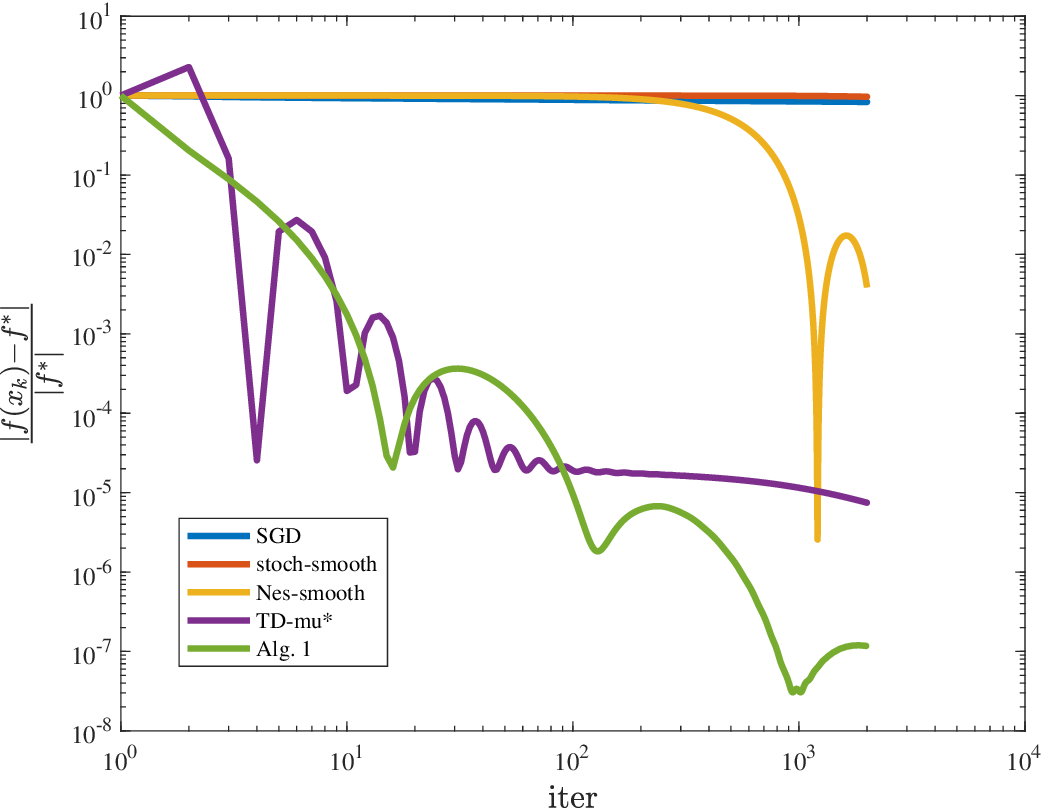}}
\hspace{8mm}
\subfloat[$\mathcal{R}(y) = \|y\|_1$ and $\eta = 10$.]{\label{fig5:sfig2}\includegraphics[scale=0.32]{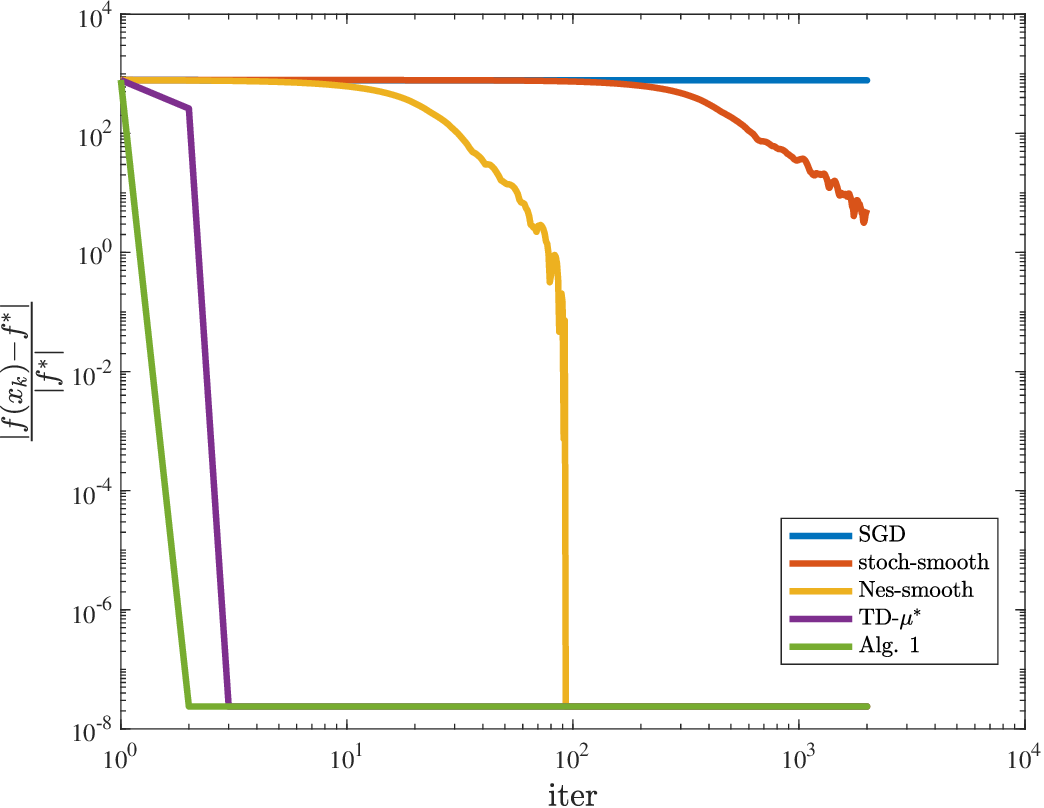}}
\caption{Numerical results for \textit{MaxCut} problem with different choice of $\mathcal{R}(y)$ in Subsection~\ref{sec:maxcut}.}
\label{sdp-appendix}
\end{figure}

\textbullet \, Figure \ref{sdp-appendix} illustrates the simulation results for the minimizing problem \eqref{maxcut dual with penalty} different choice of $\eta$.

\begin{table*}[t]
\centering
\caption{Applications of Selected Regression Datasets from LIBSVM.}
\small
\begin{tabular}{p{2.7cm}|p{9cm}|p{4cm}}
\hline
\textbf{Dataset} & \textbf{Application / Context} & \textbf{Reference} \\
\hline
\textbf{"mg" (Mackey-Glass)} {(Fig. \ref{ADMM1}--\ref{ADMM2})} & Time-series regression problem modeling a gas furnace system. The goal is to predict CO\textsubscript{2} concentration based on methane input and past outputs. Used in system identification and control engineering. & LIBSVM dataset; Box–Jenkins studies \cite{chang2011libsvm} \\
\hline
\textbf{"abalone"} {(Fig. \ref{ADMM1}--\ref{ADMM2})} & Predict the age of abalone (a type of mollusk) from physical measurements such as shell length, diameter, and weight. Used in marine biology for studying growth and population dynamics. & UCI Abalone dataset; LIBSVM \cite{chang2011libsvm} \\
\hline
\textbf{"cpusmall"} & Predict the relative performance of CPUs from their hardware specifications (e.g., memory, clock speed). Used in benchmarking and performance modeling in computer architecture. & StatLib CPU performance dataset; LIBSVM \cite{chang2011libsvm} \\
\hline
\textbf{"triazines"} & Predict the biological activity of triazine chemical compounds using molecular descriptors. Commonly used in chemoinformatics and quantitative structure activity relationship (QSAR) modeling for drug discovery. & QSAR modeling studies; Statlog; LIBSVM \cite{chang2011libsvm} \\
\hline
\end{tabular}
\label{table2}
\end{table*}

\begin{figure*}[h!] 
\centering
\begin{minipage}{0.32\textwidth}
\includegraphics[width=\linewidth]{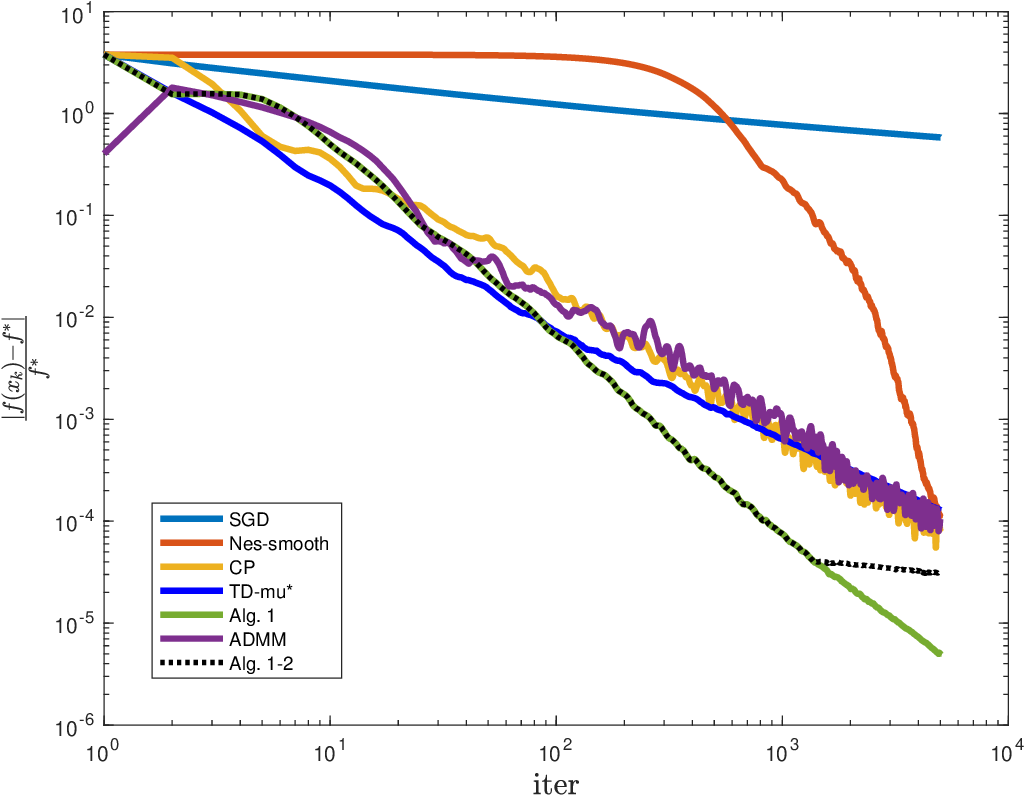}
        \subcaption{setting similar to Fig. \ref{fig1:sfig1}.}
  \label{admm1}
\end{minipage}%
\hspace{-1mm}
\begin{minipage}{0.32\textwidth}
\includegraphics[width=\linewidth]{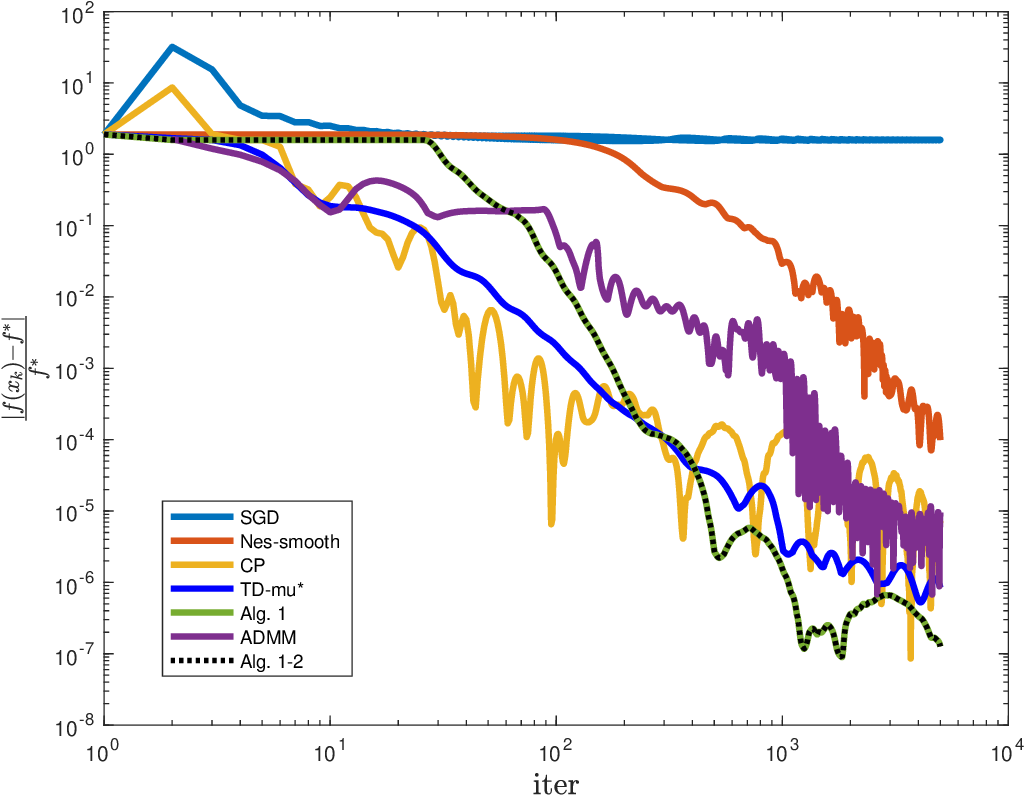}
        \subcaption{real dataset (``abalone").}
  \label{admm2}
\end{minipage}%
\hspace{-1mm}
\begin{minipage}{0.32\textwidth}
\includegraphics[width=\linewidth]{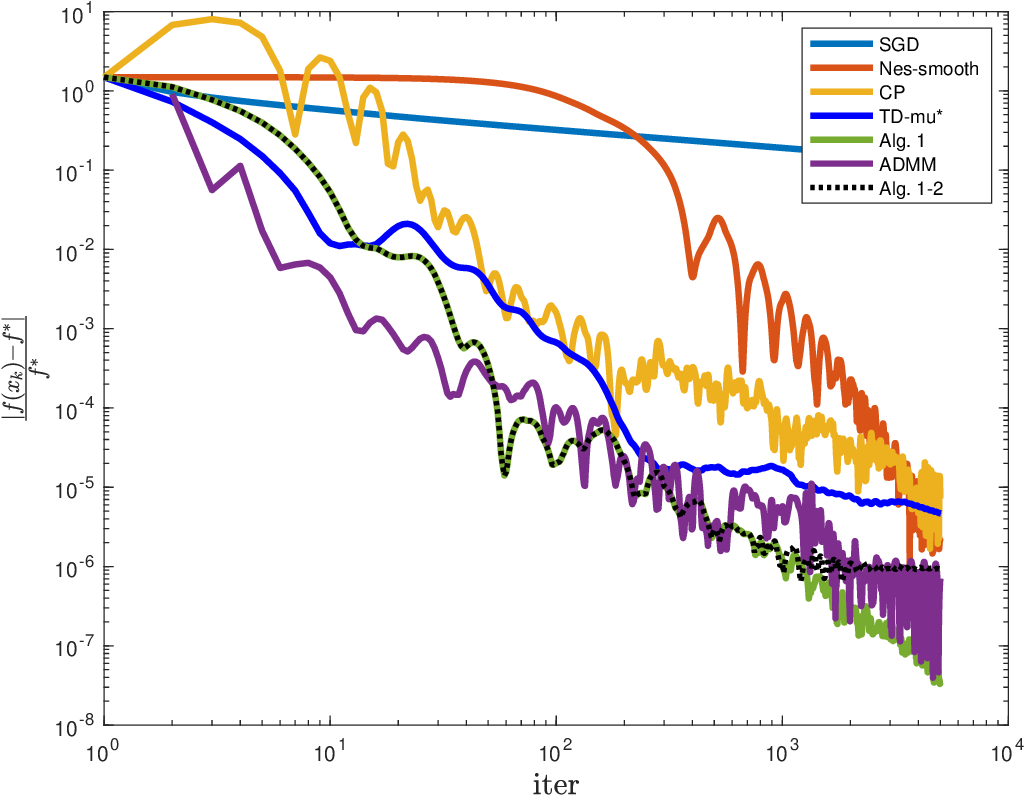}
        \subcaption{real dataset (``mg").}
  \label{admm3}
\end{minipage}%
\caption{Numerical results for $\ell_1$-$\ell_1$ regression problem as in \eqref{regression problem}.}
\label{ADMM1}
\end{figure*}

\begin{figure*}[h!] 
\centering
\begin{minipage}{0.32\textwidth}
\includegraphics[width=\linewidth]{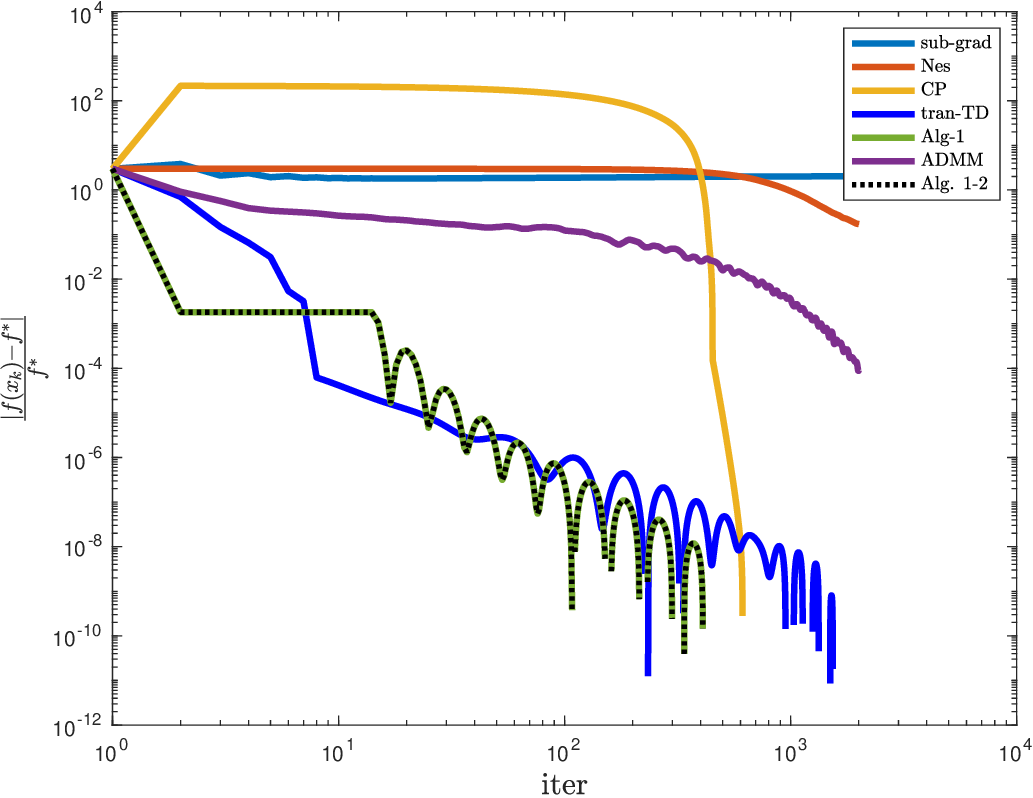}
        \subcaption{setting similar to Fig. \ref{fig1:sfig2}.}
  \label{admm4}
\end{minipage}%
\hspace{-1mm}
\begin{minipage}{0.32\textwidth}
\includegraphics[width=\linewidth]{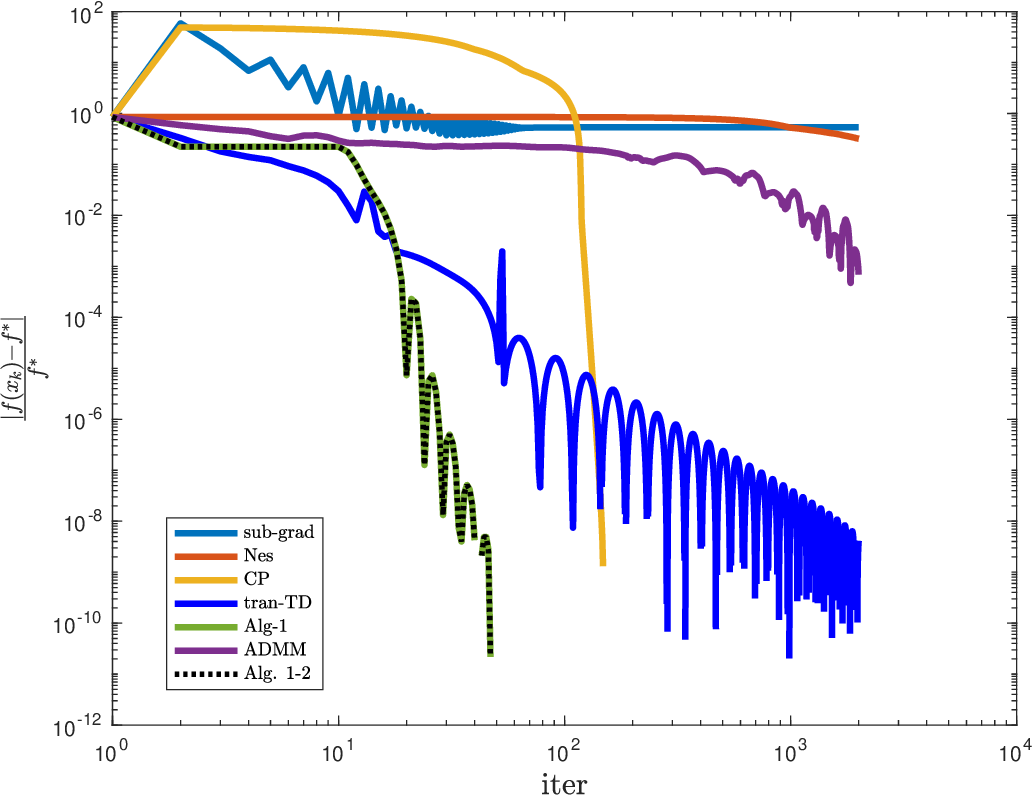}
        \subcaption{real dataset (``abalone").}
  \label{admm5}
\end{minipage}%
\hspace{-1mm}
\begin{minipage}{0.32\textwidth}
\includegraphics[width=\linewidth]{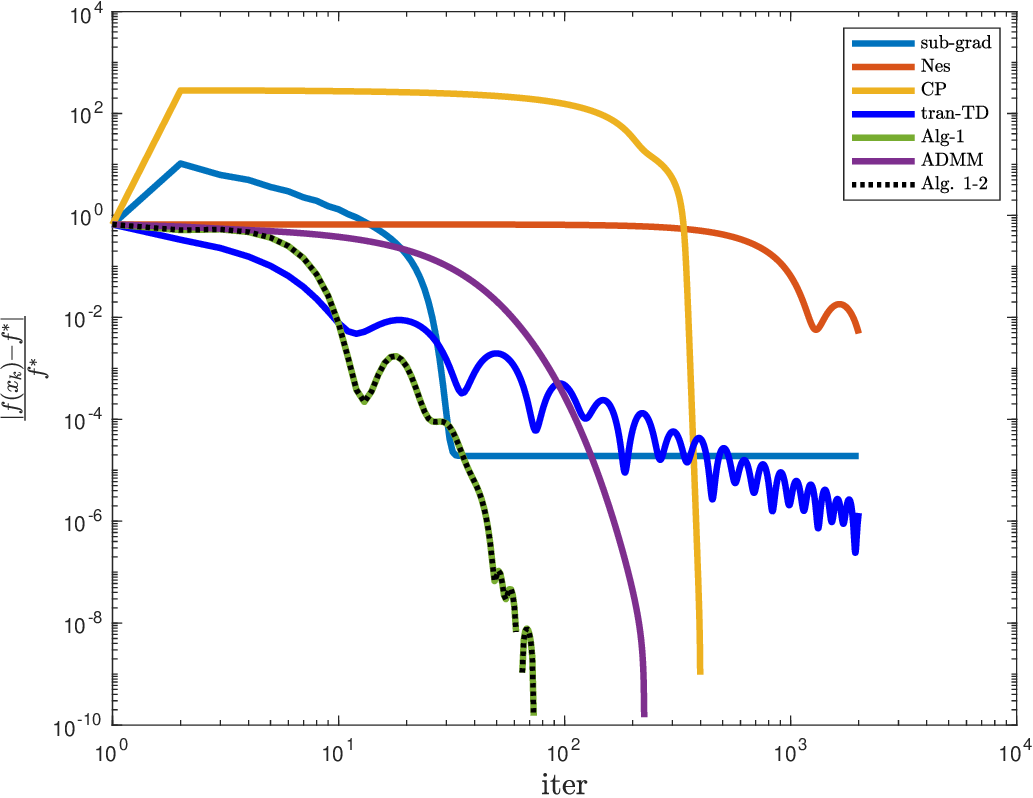}
        \subcaption{real dataset (``mg").}
  \label{admm6}
\end{minipage}%
\caption{Numerical results for $\ell_1$-$\ell_2$ regression problem as in \eqref{regression problem}.}
\label{ADMM2}
\end{figure*}

\textbullet \,\  Compared with ADMM: it seems that \eqref{main problem} is suitable for ADMM algorithm. Then, in Figures \ref{ADMM1} and \ref{ADMM2} we have also implemented our proposed algorithm (both versions with $c = 0$ and $c = {\epsilon}/{L_{f}^2}$ in the smoothing rule) together with the current four existing approaches in the paper and ADMM from Algorithm 2 in \cite{beck2017first} on two real datasets from the regression data of ``abalone'' and ``mg'' from LIBSVM \cite{chang2011libsvm}. We run the ADMM with three different penalty parameters ($\varrho = 0.1,1,10$) and report the best one.

\begin{figure}[h!]
\centering
\subfloat[$\ell_1$-$\ell_1$-regularized LASSO.]{\label{fig11:sfig1}\includegraphics[scale=0.32]{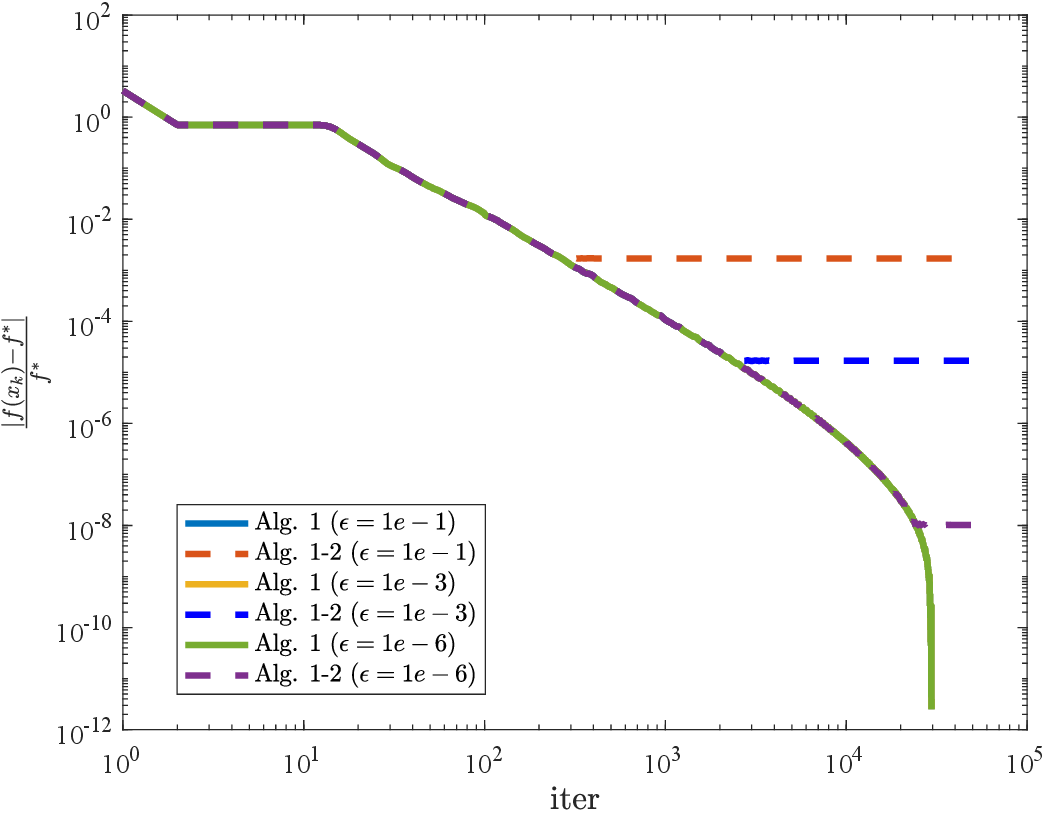}}
\hspace{8mm}
\subfloat[$\ell_2$-$\ell_1$-regularized LASSO.]{\label{fig11:sfig2}\includegraphics[scale=0.32]{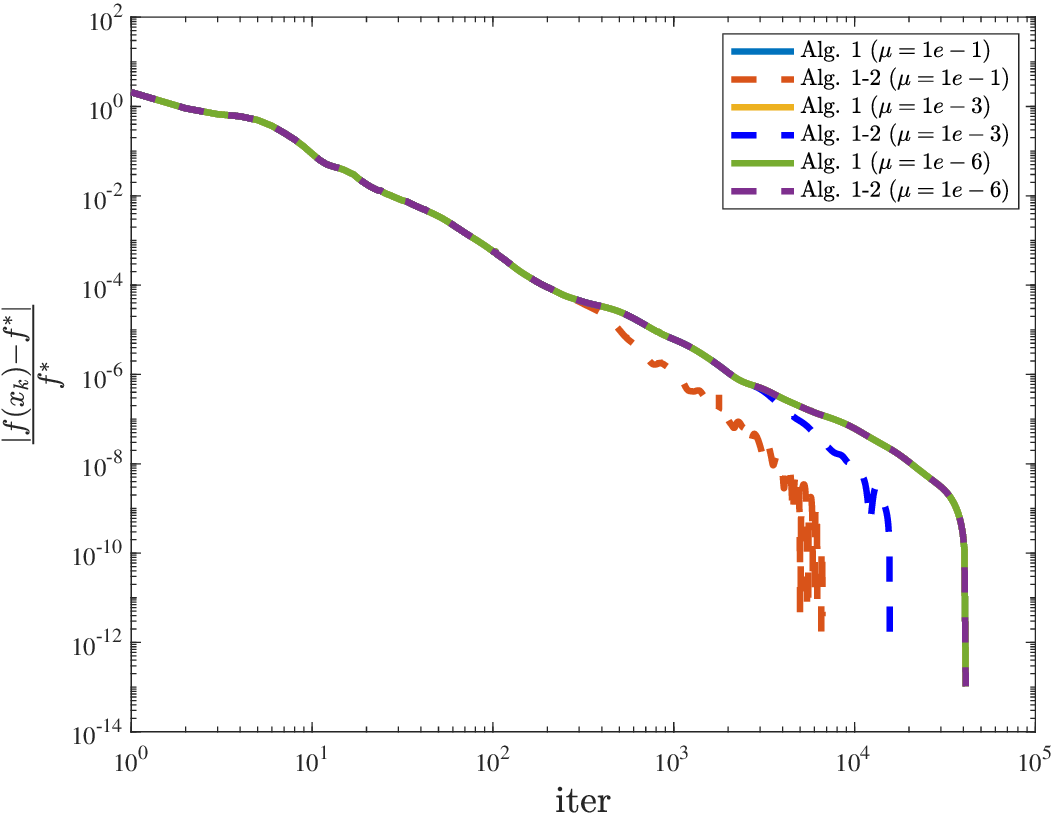}}
\caption{Numerical results for regression problems with different values of $\varepsilon$.}
\label{LASSO-compare-epsilon}
\end{figure}

\textbullet \,\ Figure \ref{LASSO-compare-epsilon} contains plots for the regression problem \eqref{regression problem} with $i = 1,2$ and correlated data (setting similar to Figures~\ref{fig1:sfig1} and \ref{fig1:sfig2}). All algorithms are initialized at the same randomly chosen point. This includes comparisons of Algorithm \ref{alg:Algorithm1}  with constant term $c = {\epsilon}/{L_{f}^2}$ (Theorem~\ref{thm:global convegence}) and various values of $\epsilon$, as well as the case where the constant term $c = 0$.

\textbullet \,\ Real-world applications and dimentionality: Indeed, the applicability of the proposed method is an important consideration. The widespread use of nonsmooth convex optimization across various domains motivated us to include 16 diverse numerical examples in the paper. Relevant applications in image processing, image reconstruction, system identification, and control, particularly those involving the summation of nonsmooth terms, are discussed with references in Numerical results section. Additionally, in the Additional Numerical Results, we include simulations using real-world datasets "mg" and "abalone" from the LIBSVM regression collection. To provide further comparison, we perform 14 additional simulations by increasing the data dimensionality (5-20 times larger) under similar settings across various problem types where we report the optimality gap $\frac{|f(x_k) - f(x^*)|}{|f(x^*)|}$ at 10 equally spaced points throughout iterations. We also included new real-world datasets, "cpusmall" and "trianzime," from LIBSVM. Table \ref{table2} summarizes the application domain of each dataset. 
\newpage \onecolumn
\textbullet \,\ \textbf{$\ell_1$-$\ell_1$-regularized LASSO}\

\noindent $n = 1000, \, m = 100$
\begin{table*}[h!]
\begin{adjustbox}{margin=-40pt 0pt 0pt 0mm} 
\scriptsize
\begin{tabular}{l|l|l|l|l|l|l|l|l|l|l}
\multicolumn{1}{l|}{Algorithm}          & point \#1 & point \#2 & point \#3 & point \#4 & point \#5 & point \#6 & point \#7 & point \#8 & point \#9 & point \#10 \\ \cline{1-1}
\hline
SGD                                      & \bf 7.74 & 2.27$\cdot10^{-2}$ & 1.69$\cdot10^{-2}$ & 1.45$\cdot10^{-2}$ & 1.28$\cdot10^{-2}$ & 1.19$\cdot10^{-2}$ & 1.13$\cdot10^{-2}$ & 1.01$\cdot10^{-2}$ & 9.76$\cdot10^{-3}$ & 9.36$\cdot10^{-3}$ \\
\hline
Nes-smooth                               & \bf7.74 & 3.43 & 6.47$\cdot10^{-1}$ & 2.79$\cdot10^{-1}$ & 1.44$\cdot10^{-1}$ & 7.50$\cdot10^{-2}$ & 3.89$\cdot10^{-2}$ & 1.96$\cdot10^{-2}$ & 1.01$\cdot10^{-2}$ & 4.24$\cdot10^{-3}$ \\
\hline
CP                                       & \bf 7.74 & 3.84$\cdot10^{-4}$ & 1.39$\cdot10^{-4}$ & 1.27$\cdot10^{-4}$ & 4.85$\cdot10^{-5}$ & 5.46$\cdot10^{-5}$ & 7.41$\cdot10^{-5}$ & 2.55$\cdot10^{-5}$ & 4.44$\cdot10^{-5}$ & 1.58$\cdot10^{-5}$ \\
\hline
TD-$\mu^*$                               & \bf 7.74 & 2.26$\cdot10^{-4}$ & 4.24$\cdot10^{-5}$ & 2.26$\cdot10^{-5}$ & 1.56$\cdot10^{-5}$ & 1.27$\cdot10^{-5}$ & 1.01$\cdot10^{-5}$ & 8.94$\cdot10^{-6}$ & 7.40$\cdot10^{-6}$ & 6.72$\cdot10^{-6}$ \\
\hline
Alg. 1                                   & \bf 7.74 & \bf 7.37 $\bf \cdot10^{-5}$ & \bf 1.55$\bf \cdot10^{-5}$ & \bf 5.04$\bf \cdot10^{-6}$ & \bf 2.35$\bf \cdot10^{-6}$ & \bf 1.49$\bf \cdot10^{-6}$ & \bf 1.01$\bf \cdot10^{-6}$ & \bf 6.94$\bf \cdot10^{-7}$ & \bf 5.10$\bf \cdot10^{-7}$ & \bf 3.53$\bf \cdot10^{-7}$ \\
\end{tabular}
\end{adjustbox}
\end{table*}

\noindent $n = 100, \, m = 1000$
\begin{table*}[h!]
\begin{adjustbox}{margin=-40pt 0pt 0pt 0mm} 
\scriptsize
\begin{tabular}{l|l|l|l|l|l|l|l|l|l|l}
 \multicolumn{1}{l|}{Algorithm} & point \#1 & point \#2 & point \#3 & point \#4 & point \#5 & point \#6 & point \#7 & point \#8 & point \#9 & point \#10 \\ 
 \hline
 SGD            & \bf 17.79 & 1.37$\cdot10^{-3}$ & 8.91$\cdot10^{-4}$ & 4.82$\cdot10^{-4}$ & 3.57$\cdot10^{-4}$ & 3.43$\cdot10^{-4}$ & 2.60$\cdot10^{-4}$ & 2.35$\cdot10^{-4}$ & 2.01$\cdot10^{-4}$ & 2.15$\cdot10^{-4}$ \\
 \hline
 Nes-smooth     & \bf 17.79 & 1.75 & 3.05$\cdot10^{-1}$ & 2.40$\cdot10^{-2}$ & \bf 3.88$\bf \cdot10^{-5}$ & 3.02$\cdot10^{-5}$ & 3.00$\cdot10^{-5}$ & 3.01$\cdot10^{-5}$ & 3.00$\cdot10^{-5}$ & 3.01$\cdot10^{-5}$ \\
 \hline
  CP             & \bf 17.79 & 1.56$\cdot10^{-3}$ & 1.69$\cdot10^{-3}$ & 1.07$\cdot10^{-3}$ & 8.38$\cdot10^{-4}$ & 6.70$\cdot10^{-4}$ & 6.77$\cdot10^{-4}$ & 5.31$\cdot10^{-4}$ & 4.49$\cdot10^{-4}$ & 3.78$\cdot10^{-4}$ \\
 \hline
 TD-$\mu^*$     & \bf 17.79 & 1.17$\cdot10^{-3}$ & 6.65$\cdot10^{-4}$ & 5.02$\cdot10^{-4}$ & 3.90$\cdot10^{-4}$ & 3.13$\cdot10^{-4}$ & 2.54$\cdot10^{-4}$ & 2.14$\cdot10^{-4}$ & 1.91$\cdot10^{-4}$ & 1.74$\cdot10^{-4}$ \\
 \hline
 Alg. 1        & \bf 17.79 & \bf 5.90$\bf \cdot10^{-4}$ & \bf 1.61$\bf \cdot10^{-4}$ & \bf 7.64$\bf \cdot10^{-5}$ & 4.26$\cdot10^{-5}$ & \bf 2.73$\bf \cdot10^{-5}$ & \bf 1.89$\bf \cdot10^{-5}$ & \bf 1.39$\bf \cdot10^{-5}$ & \bf 1.07$\bf \cdot10^{-5}$ & \bf 8.52$\bf \cdot10^{-6}$ \\
\end{tabular}
\end{adjustbox}
\end{table*}

\noindent $n = 2000, \, m = 500$
\begin{table*}[h!]
\begin{adjustbox}{margin=-40pt 0pt 0pt 0mm} 
\scriptsize
\begin{tabular}{l|l|l|l|l|l|l|l|l|l|l}
 \multicolumn{1}{l|}{Algorithm} & point \#1 & point \#2 & point \#3 & point \#4 & point \#5 & point \#6 & point \#7 & point \#8 & point \#9 & point \#10 \\ \cline{2-2}
  \hline
 SGD            & \bf 5.28 & 1.63$\cdot10^{-2}$ & 1.19$\cdot10^{-2}$ & 9.75$\cdot10^{-3}$ & 8.45$\cdot10^{-3}$ & 7.67$\cdot10^{-3}$ & 7.07$\cdot10^{-3}$ & 6.45$\cdot10^{-3}$ & 5.85$\cdot10^{-3}$ & 5.48$\cdot10^{-3}$ \\
  \hline
 Nes-smooth     & \bf 5.28 & 2.46 & 4.80$\cdot10^{-1}$ & 2.11$\cdot10^{-1}$ & 1.05$\cdot10^{-1}$ & 5.46$\cdot10^{-2}$ & 2.83$\cdot10^{-2}$ & 1.53$\cdot10^{-2}$ & 7.30$\cdot10^{-3}$ & 3.41$\cdot10^{-3}$ \\
  \hline
 CP             & \bf 5.28 & 1.22$\cdot10^{-4}$ & 4.21$\cdot10^{-5}$ & 2.43$\cdot10^{-5}$ & 1.31$\cdot10^{-5}$ & 8.99$\cdot10^{-6}$ & 7.86$\cdot10^{-6}$ & 4.77$\cdot10^{-6}$ & 3.40$\cdot10^{-6}$ & 2.92$\cdot10^{-6}$ \\
  \hline
 TD-$\mu^*$     & \bf 5.28 & 6.05$\cdot10^{-5}$ & 2.62$\cdot10^{-5}$ & 1.66$\cdot10^{-5}$ & 1.19$\cdot10^{-5}$ & 9.45$\cdot10^{-6}$ & 7.92$\cdot10^{-6}$ & 6.95$\cdot10^{-6}$ & 5.92$\cdot10^{-6}$ & 5.29$\cdot10^{-6}$ \\
  \hline
 Alg. 1        & \bf 5.28 & \bf 2.30$\bf \cdot10^{-5}$ & \bf 4.83$\bf \cdot10^{-6}$ & \bf 2.35$\bf \cdot10^{-6}$ & \bf 1.39$\bf \cdot10^{-6}$ & \bf 1.02$\bf \cdot10^{-6}$ & \bf 7.62$\bf \cdot10^{-7}$ & \bf 5.83$\bf \cdot10^{-7}$ & \bf 4.72$\bf \cdot10^{-7}$ & \bf 3.93$\bf \cdot10^{-7}$ \\
\end{tabular}
\end{adjustbox}
\end{table*}

\noindent $n = 500, \, m = 2000$
\begin{table*}[h!]
\begin{adjustbox}{margin=-40pt 0pt 0pt 0mm} 
\scriptsize
\begin{tabular}{l|l|l|l|l|l|l|l|l|l|l}
 \multicolumn{1}{l|}{Algorithm} & point \#1 & point \#2 & point \#3 & point \#4 & point \#5 & point \#6 & point \#7 & point \#8 & point \#9 & point \#10 \\ \cline{2-2}
 \hline
 SGD            & \bf 2.29 & \bf 4.29$\bf \cdot10^{-5}$ & \bf 1.58$\bf \cdot10^{-5}$ & 9.28$\cdot10^{-6}$ & 6.58$\cdot10^{-6}$ & 5.70$\cdot10^{-6}$ & 4.64$\cdot10^{-6}$ & 4.01$\cdot10^{-6}$ & 3.60$\cdot10^{-6}$ & 3.73$\cdot10^{-6}$ \\
 \hline
 Nes-smooth     & \bf 2.29 & 5.63$\cdot10^{-1}$ & 9.02$\cdot10^{-2}$ & 3.23$\cdot10^{-2}$ & 1.51$\cdot10^{-2}$ & 6.63$\cdot10^{-3}$ & 1.89$\cdot10^{-3}$ & 1.44$\cdot10^{-4}$ & 6.14$\cdot10^{-6}$ & 4.24$\cdot10^{-6}$ \\
 \hline
 CP             & \bf 2.29 & 4.77$\cdot10^{-4}$ & 2.50$\cdot10^{-4}$ & 1.87$\cdot10^{-4}$ & 1.49$\cdot10^{-4}$ & 1.17$\cdot10^{-4}$ & 9.05$\cdot10^{-5}$ & 7.90$\cdot10^{-5}$ & 6.35$\cdot10^{-5}$ & 6.03$\cdot10^{-5}$ \\
 \hline
 TD-$\mu^*$     & \bf 2.29 & 4.37$\cdot10^{-4}$ & 1.84$\cdot10^{-4}$ & 1.21$\cdot10^{-4}$ & 9.20$\cdot10^{-5}$ & 7.50$\cdot10^{-5}$ & 6.34$\cdot10^{-5}$ & 5.46$\cdot10^{-5}$ & 4.77$\cdot10^{-5}$ & 4.23$\cdot10^{-5}$ \\
 \hline
 Alg. 1         & \bf 2.29 & 8.54$\cdot10^{-5}$ & 2.06$\cdot10^{-5}$ & \bf 8.55$\bf \cdot10^{-6}$ & \bf 4.58$\bf \cdot10^{-6}$ & \bf 2.77$\bf \cdot10^{-6}$ & \bf 1.82$\bf \cdot10^{-6}$ & \bf 1.27$\bf \cdot10^{-6}$ & \bf 9.21$\bf \cdot10^{-7}$ & \bf 6.86$\bf \cdot10^{-7}$ \\
\end{tabular}
\end{adjustbox}
\end{table*}

cpusmall dataset $n = 12,\, m = 8192$
\begin{table*}[h!]
\begin{adjustbox}{margin=-40pt 0pt 0pt 0mm} 
\scriptsize
\begin{tabular}{l|l|l|l|l|l|l|l|l|l|l}
 \multicolumn{1}{l|}{Algorithm} & point \#1 & point \#2 & point \#3 & point \#4 & point \#5 & point \#6 & point \#7 & point \#8 & point \#9 & point \#10 \\ \cline{2-2}
  \hline
 SGD            & \bf 1.62$\bf \cdot10^5$ & 3.56$\cdot10^2$ & 7.61$\cdot10^1$ & 7.54$\cdot10^1$ & 3.46$\cdot10^2$ & 3.44$\cdot10^2$ & 7.32$\cdot10^1$ & 7.24$\cdot10^1$ & 3.36$\cdot10^2$ & 3.33$\cdot10^2$ \\
  \hline
 Nes-smooth     & \bf 1.62$\bf \cdot10^5$ & 1.62$\cdot10^5$ & 1.62$\cdot10^5$ & 1.62$\cdot10^5$ & 1.62$\cdot10^5$ & 1.62$\cdot10^5$ & 1.62$\cdot10^5$ & 1.62$\cdot10^5$ & 1.62$\cdot10^5$ & 1.62$\cdot10^5$ \\
  \hline
 CP             & \bf 1.62$\bf \cdot10^5$ & 1.64$\cdot10^3$ & 7.07$\cdot10^2$ & 3.43$\cdot10^2$ & 4.45$\cdot10^2$ & 7.03$\cdot10^2$ & 4.57$\cdot10^2$ & 2.86$\cdot10^2$ & 6.77$\cdot10^2$ & 6.48$\cdot10^2$ \\
  \hline
 TD-$\mu^*$     & \bf 1.62$\bf \cdot10^5$ & 6.99$\cdot10^1$ & 4.43$\cdot10^1$ & 1.62$\cdot10^1$ & 1.36$\cdot10^1$ & 1.86$\cdot10^1$ & 1.30$\cdot10^1$ & 9.40$\cdot10^0$ & 1.41$\cdot10^1$ & 1.36$\cdot10^1$ \\
  \hline
 Alg. 1         & \bf 1.62$\bf \cdot10^5$ & \bf 1.43$\bf \cdot10^0$ & \bf 1.14$\bf \cdot10^0$ & \bf 9.00$\bf \cdot10^{-1}$ & \bf 7.27$\bf \cdot10^{-1}$ & \bf 5.97$\bf \cdot10^{-1}$ & \bf 4.98$\bf \cdot10^{-1}$ & \bf 4.20$\bf \cdot10^{-1}$ & \bf 3.60$\bf \cdot10^{-1}$ & \bf 3.13$\bf \cdot10^{-1}$ \\
\end{tabular}
\end{adjustbox}
\end{table*}

trianzines dataset $n = 60,\, m = 186$
\begin{table*}[h!]
\begin{adjustbox}{margin=-40pt 0pt 0pt 0mm} 
\scriptsize
\begin{tabular}{l|l|l|l|l|l|l|l|l|l|l}
 \multicolumn{1}{l|}{Algorithm} & point \#1 & point \#2 & point \#3 & point \#4 & point \#5 & point \#6 & point \#7 & point \#8 & point \#9 & point \#10 \\ 
  \hline
 SGD            & \bf 2.16$\bf \cdot10^2$ & 5.02$\cdot10^{-1}$ & 7.84$\cdot10^{-2}$ & 9.15$\cdot10^{-2}$ & 6.76$\cdot10^{-2}$ & 8.46$\cdot10^{-2}$ & 4.20$\cdot10^{-2}$ & 4.92$\cdot10^{-2}$ & 3.91$\cdot10^{-2}$ & 2.60$\cdot10^{-2}$ \\
  \hline
 Nes-smooth     & \bf 2.16$\bf \cdot10^2$ & 1.25$\cdot10^1$ & 5.22$\cdot10^0$ & 9.11$\cdot10^{-1}$ & 5.03$\cdot10^{-1}$ & 2.80$\cdot10^{-1}$ & 1.74$\cdot10^{-1}$ & 1.19$\cdot10^{-1}$ & 7.86$\cdot10^{-2}$ & 4.46$\cdot10^{-2}$ \\
  \hline
 CP             & \bf 2.16$\bf \cdot10^2$ & 1.03$\cdot10^{-1}$ & 4.34$\cdot10^{-2}$ & 2.97$\cdot10^{-2}$ & 2.61$\cdot10^{-2}$ & 1.55$\cdot10^{-2}$ & 1.54$\cdot10^{-2}$ & 1.44$\cdot10^{-2}$ & 8.14$\cdot10^{-3}$ & 9.77$\cdot10^{-3}$ \\
  \hline
 TD-$\mu^*$     & \bf 2.16$\bf \cdot10^2$ & 2.21$\cdot10^{-2}$ & 9.72$\cdot10^{-3}$ & 4.57$\cdot10^{-3}$ & 4.05$\cdot10^{-3}$ & 3.36$\cdot10^{-3}$ & 2.55$\cdot10^{-3}$ & 1.99$\cdot10^{-3}$ & 1.87$\cdot10^{-3}$ & 1.59$\cdot10^{-3}$ \\
  \hline
 Alg. 1         & \bf 2.16$\bf \cdot10^2$ & \bf 1.71$\bf \cdot10^{-2}$ & \bf 4.48$\bf \cdot10^{-3}$ & \bf 1.72$\bf \cdot10^{-3}$ & \bf 8.98$\bf \cdot10^{-4}$ & \bf 6.66$\bf \cdot10^{-4}$ & \bf 4.08$\bf \cdot10^{-4}$ & \bf 3.58$\bf \cdot10^{-4}$ & \bf 2.63$\bf \cdot10^{-4}$ & \bf 1.86$\bf \cdot10^{-4}$ \\
\end{tabular}
\end{adjustbox}
\end{table*}
\vspace{10cm}

\textbullet \,\ \textbf{$\ell_2$-$\ell_1$-regularized LASSO}\

\noindent $n = 1000, \, m = 100$
\begin{table*}[h!]
\begin{adjustbox}{margin=-40pt 0pt 0pt 0mm} 
\scriptsize
\begin{tabular}{l|l|l|l|l|l|l|l|l|l|l}
 \multicolumn{1}{l|}{Algorithm} & point \#1 & point \#2 & point \#3 & point \#4 & point \#5 & point \#6 & point \#7 & point \#8 & point \#9 & point \#10 \\ 
 \hline
 SGD            & \bf 7.57 & 5.75$\cdot 10^{-2}$ & 4.42$\cdot 10^{-2}$ & 3.89$\cdot 10^{-2}$ & 3.58$\cdot 10^{-2}$ & 3.38$\cdot 10^{-2}$ & 3.24$\cdot 10^{-2}$ & 3.12$\cdot 10^{-2}$ & 3.03$\cdot 10^{-2}$ & 2.94$\cdot 10^{-2}$ \\
 \hline
 Nes-smooth     & \bf 7.57 & 7.20$\cdot 10^{0}$  & 6.19$\cdot 10^{0}$  & 4.76$\cdot 10^{0}$  & 3.32$\cdot 10^{0}$  & 2.39$\cdot 10^{0}$  & 1.66$\cdot 10^{0}$  & 1.15$\cdot 10^{0}$  & 8.18$\cdot 10^{-1}$ & 6.66$\cdot 10^{-1}$ \\
 \hline
 CP             & \bf 7.57 & 8.60$\cdot 10^{0}$  & 5.86$\cdot 10^{-2}$ & 1.98$\cdot 10^{-2}$ & 9.62$\cdot 10^{-3}$ & 5.73$\cdot 10^{-3}$ & 3.58$\cdot 10^{-3}$ & 2.51$\cdot 10^{-3}$ & 1.76$\cdot 10^{-3}$ & 1.27$\cdot 10^{-3}$ \\
 \hline
 TD-$\mu^*$     & \bf 7.57 & \bf 2.76$\bf \cdot 10^{-3}$ & 1.37$\cdot 10^{-3}$ & 9.06$\cdot 10^{-4}$ & 6.74$\cdot 10^{-4}$ & 5.36$\cdot 10^{-4}$ & 4.45$\cdot 10^{-4}$ & 3.82$\cdot 10^{-4}$ & 3.34$\cdot 10^{-4}$ & 2.97$\cdot 10^{-4}$ \\
 \hline
 Alg. 1         & \bf 7.57 & 3.63$\cdot 10^{-3}$ & \bf 9.11$\bf \cdot 10^{-4}$ & \bf 4.06$\bf \cdot 10^{-4}$ & \bf 2.28$\bf \cdot 10^{-4}$ & \bf 1.46$\bf \cdot 10^{-4}$ & \bf 1.01$\bf \cdot 10^{-4}$ & \bf 7.43$\bf \cdot 10^{-5}$ & \bf 5.70$\bf \cdot 10^{-5}$ & \bf 4.50$\bf \cdot 10^{-5}$ \\
\end{tabular}
\end{adjustbox}
\end{table*}

\noindent $n = 100, \, m = 1000$
\begin{table*}[h!]
\begin{adjustbox}{margin=-40pt 0pt 0pt 0mm} 
\scriptsize
\begin{tabular}{l|l|l|l|l|l|l|l|l|l|l}
 \multicolumn{1}{l|}{Algorithm} & point \#1 & point \#2 & point \#3 & point \#4 & point \#5 & point \#6 & point \#7 & point \#8 & point \#9 & point \#10 \\ 
 \hline
 SGD            & \bf 1.91$\bf \cdot10^2$ & 1.61$\cdot10^{-2}$ & \bf 8.33$\bf \cdot10^{-9}$ & \bf 8.33$\bf \cdot10^{-9}$ & \bf 8.33$\bf \cdot10^{-9}$ & \bf 8.33$\bf \cdot10^{-9}$ & \bf 8.33$\bf \cdot10^{-9}$ & \bf 8.33$\bf \cdot10^{-9}$ & \bf 8.33$\bf \cdot10^{-9}$ & \bf 8.33$\bf \cdot10^{-9}$ \\
 \hline
 Nes-smooth     & \bf 1.91$\bf \cdot10^2$ & 1.79$\cdot10^2$ & 1.45$\cdot10^2$ & 9.21$\cdot10^1$ & 3.15$\cdot10^1$ & 2.10$\cdot10^1$ & 1.82$\cdot10^1$ & 9.64$\cdot10^0$ & 9.05$\cdot10^0$ & 7.50$\cdot10^0$ \\
 \hline
 CP             & \bf 1.91$\bf \cdot10^2$ & 2.39$\cdot10^3$ & 1.06$\cdot10^0$ & \bf 8.33$\bf \cdot10^{-9}$ & \bf 8.33$\bf \cdot10^{-9}$ & \bf 8.33$\bf \cdot10^{-9}$ & \bf 8.33$\bf \cdot10^{-9}$ & \bf 8.33$\bf \cdot10^{-9}$ & \bf 8.33$\bf \cdot10^{-9}$ & \bf 8.33$\bf \cdot10^{-9}$ \\
 \hline
 TD-$\mu^*$     & \bf 1.91$\bf \cdot10^2$ & \bf 8.33$\bf \cdot10^{-9}$ & \bf 8.33$\bf \cdot10^{-9}$ & \bf 8.33$\bf \cdot10^{-9}$ & \bf 8.33$\bf \cdot10^{-9}$ & \bf 8.33$\bf \cdot10^{-9}$ & \bf 8.33$\bf \cdot10^{-9}$ & \bf 8.33$\bf \cdot10^{-9}$ & \bf 8.33$\bf \cdot10^{-9}$ & \bf 8.33$\bf \cdot10^{-9}$ \\
 \hline
 Alg. 1         & \bf 1.91$\bf \cdot10^2$ & \bf 8.33$\bf \cdot10^{-9}$ & \bf 8.33$\bf \cdot10^{-9}$ & \bf 8.33$\bf \cdot10^{-9}$ & \bf 8.33$\bf \cdot10^{-9}$ & \bf 8.33$\bf \cdot10^{-9}$ & \bf 8.33$\bf \cdot10^{-9}$ & \bf 8.33$\bf \cdot10^{-9}$ & \bf 8.33$\bf \cdot10^{-9}$ & \bf 8.33$\bf \cdot10^{-9}$ \\
\end{tabular}
\end{adjustbox}
\end{table*}

\noindent $n = 2000, \, m = 500$
\begin{table*}[h!]
\begin{adjustbox}{margin=-40pt 0pt 0pt 0mm} 
\scriptsize
\begin{tabular}{l|l|l|l|l|l|l|l|l|l|l}
 \multicolumn{1}{l|}{Algorithm} & point \#1 & point \#2 & point \#3 & point \#4 & point \#5 & point \#6 & point \#7 & point \#8 & point \#9 & point \#10 \\ 
  \hline
 SGD            & \bf 4.04 & 2.92$\cdot 10^{-2}$ & 2.30$\cdot 10^{-2}$ & 2.01$\cdot 10^{-2}$ & 1.86$\cdot 10^{-2}$ & 1.76$\cdot 10^{-2}$ & 1.67$\cdot 10^{-2}$ & 1.61$\cdot 10^{-2}$ & 1.56$\cdot 10^{-2}$ & 1.52$\cdot 10^{-2}$ \\
  \hline
 Nes-smooth     & \bf 4.04  & 3.96$\cdot 10^{0}$  & 3.72$\cdot 10^{0}$  & 3.33$\cdot 10^{0}$  & 2.84$\cdot 10^{0}$  & 2.30$\cdot 10^{0}$  & 1.79$\cdot 10^{0}$  & 1.35$\cdot 10^{0}$  & 1.04$\cdot 10^{0}$  & 8.24$\cdot 10^{-1}$ \\
  \hline
 CP             & \bf 4.04  & 5.29$\cdot 10^{1}$  & 1.20$\cdot 10^{1}$  & 2.07$\cdot 10^{-1}$ & 3.14$\cdot 10^{-2}$ & 1.36$\cdot 10^{-2}$ & 6.94$\cdot 10^{-3}$ & 4.32$\cdot 10^{-3}$ & 2.94$\cdot 10^{-3}$ & 2.15$\cdot 10^{-3}$ \\
  \hline
 TD-$\mu^*$     & \bf 4.04  & \bf 1.07$\bf \cdot 10^{-3}$ & 5.23$\cdot 10^{-4}$ & 3.46$\cdot 10^{-4}$ & 2.57$\cdot 10^{-4}$ & 2.05$\cdot 10^{-4}$ & 1.70$\cdot 10^{-4}$ & 1.46$\cdot 10^{-4}$ & 1.27$\cdot 10^{-4}$ & 1.13$\cdot 10^{-4}$ \\
  \hline
 Alg. 1         & \bf 4.04  & 1.49$\cdot 10^{-3}$ & \bf 3.82$\bf \cdot 10^{-4}$ & \bf 1.71$\bf \cdot 10^{-4}$ & \bf 9.69$\bf \cdot 10^{-5}$ & \bf 6.22$\bf \cdot 10^{-5}$ & \bf 4.33$\bf \cdot 10^{-5}$ & \bf 3.22$\bf \cdot 10^{-5}$ & \bf 2.48$\bf \cdot 10^{-5}$ & \bf 1.98$\bf \cdot 10^{-5}$ \\
\end{tabular}
\end{adjustbox}
\end{table*}

\noindent $n = 500, \, m = 2000$
\begin{table*}[h!]
\begin{adjustbox}{margin=-40pt 0pt 0pt 0mm} 
\scriptsize
\begin{tabular}{l|l|l|l|l|l|l|l|l|l|l}
 \multicolumn{1}{l|}{Algorithm} & point \#1 & point \#2 & point \#3 & point \#4 & point \#5 & point \#6 & point \#7 & point \#8 & point \#9 & point \#10 \\ 
  \hline
 SGD            & \bf 2.43$\bf \cdot 10^{2}$ & 3.47$\cdot 10^{-1}$ & 9.03$\cdot 10^{-3}$ & \bf 9.27$\bf \cdot 10^{-9}$ & \bf 9.27$\bf \cdot 10^{-9}$ &\bf  9.27$\bf \cdot 10^{-9}$ & \bf 9.27$\bf \cdot 10^{-9}$ & \bf 9.27$\bf \cdot 10^{-9}$ & \bf 9.27$\bf \cdot 10^{-9}$ & \bf 9.27$\bf \cdot 10^{-9}$ \\
  \hline
 Nes-smooth     & \bf 2.43$\bf \cdot 10^{2}$ & 2.39$\cdot 10^{2}$ & 2.28$\cdot 10^{2}$ & 2.10$\cdot 10^{2}$ & 1.85$\cdot 10^{2}$ & 1.55$\cdot 10^{2}$ & 1.22$\cdot 10^{2}$ & 8.85$\cdot 10^{1}$ & 5.97$\cdot 10^{1}$ & 4.04$\cdot 10^{1}$ \\
  \hline
 CP             & \bf 2.43$\bf \cdot 10^{2}$ & 4.64$\cdot 10^{3}$ & 2.80$\cdot 10^{3}$ & 1.33$\cdot 10^{3}$ & 3.53$\cdot 10^{2}$ & 8.71$\cdot 10^{-4}$ & \bf 9.27$\bf \cdot 10^{-9}$ & \bf 9.27$\bf \cdot 10^{-9}$ & \bf 9.27$\bf \cdot 10^{-9}$ & \bf 9.27$\bf \cdot 10^{-9}$ \\
  \hline
 TD-$\mu^*$     & \bf 2.43$\bf \cdot 10^{2}$ & 8.29$\cdot 10^{-8}$ & \bf 9.27$\bf \cdot 10^{-9}$ & \bf 9.27$\bf \cdot 10^{-9}$ & \bf 9.27$\bf \cdot 10^{-9}$ & \bf 9.27$\bf \cdot 10^{-9}$ & \bf 9.27$\bf \cdot 10^{-9}$ & \bf 9.27$\bf \cdot 10^{-9}$ & \bf 9.27$\bf \cdot 10^{-9}$ & \bf 9.27$\bf \cdot 10^{-9}$ \\
  \hline
 Alg. 1         & \bf 2.43$\bf \cdot 10^{2}$ & \bf 9.25$\bf \cdot 10^{-9}$ & \bf 9.27$\bf \cdot 10^{-9}$ & \bf 9.27$\bf \cdot 10^{-9}$ & \bf 9.27$\bf \cdot 10^{-9}$ & \bf 9.27$\bf \cdot 10^{-9}$ & \bf 9.27$\bf \cdot 10^{-9}$ & \bf 9.27$\bf \cdot 10^{-9}$ &\bf  9.27$\bf \cdot 10^{-9}$ & \bf 9.27$\bf \cdot 10^{-9}$ \\
\end{tabular}
\end{adjustbox}
\end{table*}
\vspace{5cm}

cpusmall dataset $n = 12,\, m = 8192$
\begin{table*}[h!]
\begin{adjustbox}{margin=-40pt 0pt 0pt 0mm} 
\scriptsize
\begin{tabular}{l|l|l|l|l|l|l|l|l|l|l}
 \multicolumn{1}{l|}{Algorithm} & point \#1 & point \#2 & point \#3 & point \#4 & point \#5 & point \#6 & point \#7 & point \#8 & point \#9 & point \#10 \\ 
   \hline
 SGD            & \bf 9.97$\bf \cdot 10^{4}$ & 5.38$\cdot 10^{4}$ & 3.23$\cdot 10^{5}$ & 2.63$\cdot 10^{5}$ & 8.66$\cdot 10^{3}$ & 3.71$\cdot 10^{3}$ & 1.86$\cdot 10^{5}$ & 1.74$\cdot 10^{5}$ & 4.28$\cdot 10^{3}$ & 3.14$\cdot 10^{3}$ \\
   \hline
 Nes-smooth     & \bf 9.97$\bf \cdot 10^{4}$ & 9.97$\cdot 10^{4}$ & 9.97$\cdot 10^{4}$ & 9.97$\cdot 10^{4}$ & 9.97$\cdot 10^{4}$ & 9.97$\cdot 10^{4}$ & 9.97$\cdot 10^{4}$ & 9.97$\cdot 10^{4}$ & 9.97$\cdot 10^{4}$ & 9.97$\cdot 10^{4}$ \\
   \hline
 CP             & \bf 9.97$\bf \cdot 10^{4}$ & 4.97$\cdot 10^{3}$ & 2.01$\cdot 10^{3}$ & 2.66$\cdot 10^{2}$ & 9.82$\cdot 10^{2}$ & 7.48$\cdot 10^{2}$ & 2.57$\cdot 10^{2}$ & 7.99$\cdot 10^{1}$ & 1.65$\cdot 10^{2}$ & 1.13$\cdot 10^{2}$ \\
   \hline
 TD-$\mu^*$     & \bf 9.97$\bf \cdot 10^{4}$ & 4.67$\cdot 10^{1}$ & 2.20$\cdot 10^{1}$ & 1.57$\cdot 10^{1}$ & 9.63$\cdot 10^{0}$ & 5.76$\cdot 10^{0}$ & 5.82$\cdot 10^{0}$ & 7.20$\cdot 10^{0}$ & 7.53$\cdot 10^{0}$ & 6.60$\cdot 10^{0}$ \\
   \hline
 Alg. 1         & \bf 9.97$\bf \cdot 10^{4}$ & \bf 1.61$\bf \cdot 10^{-1}$ & \bf 1.31$\bf \cdot 10^{-1}$ & \bf 7.50$\bf \cdot 10^{-2}$ & \bf 5.57$\bf \cdot 10^{-2}$ & \bf 5.03$\bf \cdot 10^{-2}$ & \bf 4.84$\bf \cdot 10^{-2}$ & \bf 4.76$\bf \cdot 10^{-2}$ & \bf 4.72$\bf \cdot 10^{-2}$ & \bf 4.70$\bf \cdot 10^{-2}$ \\
\end{tabular}
\end{adjustbox}
\end{table*}

trianzines dataset $n = 60,\, m = 186$
\begin{table*}[h!]
\begin{adjustbox}{margin=-40pt 0pt 0pt 0mm} 
\scriptsize
\begin{tabular}{l|l|l|l|l|l|l|l|l|l|l}
 \multicolumn{1}{l|}{Algorithm} & point \#1 & point \#2 & point \#3 & point \#4 & point \#5 & point \#6 & point \#7 & point \#8 & point \#9 & point \#10 \\
    \hline
 SGD            & \bf 6.58$\bf \cdot 10^{1}$ & 6.57$\cdot 10^{1}$ & 6.56$\cdot 10^{1}$ & 6.56$\cdot 10^{1}$ & 6.55$\cdot 10^{1}$ & 6.55$\cdot 10^{1}$ & 6.55$\cdot 10^{1}$ & 6.55$\cdot 10^{1}$ & 6.54$\cdot 10^{1}$ & 6.54$\cdot 10^{1}$ \\
    \hline
 Nes-smooth     & \bf 6.58$\bf \cdot 10^{1}$ & 4.92$\cdot 10^{1}$ & 2.66$\cdot 10^{1}$ & 2.11$\cdot 10^{1}$ & 1.76$\cdot 10^{1}$ & 1.61$\cdot 10^{1}$ & 1.45$\cdot 10^{1}$ & 1.33$\cdot 10^{1}$ & 1.18$\cdot 10^{1}$ & 1.02$\cdot 10^{1}$ \\
    \hline
 CP             & \bf 6.58$\bf \cdot 10^{1}$ & 1.71$\cdot 10^{1}$ & 5.04$\cdot 10^{-1}$ & 2.77$\cdot 10^{-1}$ & 1.79$\cdot 10^{-1}$ & 1.16$\cdot 10^{-1}$ & 6.81$\cdot 10^{-2}$ & 4.12$\cdot 10^{-2}$ & 2.07$\cdot 10^{-2}$ & 2.85$\cdot 10^{-3}$ \\
    \hline
 TD-$\mu^*$     & \bf 6.58$\bf \cdot 10^{1}$ & 1.59$\cdot 10^{-1}$ & 5.98$\cdot 10^{-2}$ & 3.59$\cdot 10^{-2}$ & 1.88$\cdot 10^{-2}$ & 1.77$\cdot 10^{-2}$ & 1.28$\cdot 10^{-2}$ & 9.91$\cdot 10^{-3}$ & 8.18$\cdot 10^{-3}$ & 6.61$\cdot 10^{-3}$ \\
    \hline
 Alg. 1         & \bf 6.58$\bf \cdot 10^{1}$ & \bf 5.21$\bf \cdot 10^{-2}$ & \bf 3.92$\bf \cdot 10^{-4}$ & \bf 2.98$\bf \cdot 10^{-5}$ & \bf 2.61$\bf \cdot 10^{-5}$ & \bf 1.86$\bf \cdot 10^{-5}$ & \bf 9.84$\bf \cdot 10^{-6}$ & \bf 6.39$\bf \cdot 10^{-6}$ & \bf 5.92$\bf \cdot 10^{-6}$ & \bf 5.76$\bf \cdot 10^{-6}$ \\
\end{tabular}
\end{adjustbox}
\end{table*}

\textbullet \,\ \textbf{MaxCut problem} 

$\mathcal{R}(y) = \|y\|_1,\, \eta = 1$, dimension of matrix $= 300$

\begin{table*}[h!]
\begin{adjustbox}{margin=-40pt 0pt 0pt 0mm} 
\scriptsize
\begin{tabular}{l|l|l|l|l|l|l|l|l|l|l}
 \multicolumn{1}{l|}{Algorithm} & point \#1 & point \#2 & point \#3 & point \#4 & point \#5 & point \#6 & point \#7 & point \#8 & point \#9 & point \#10 \\ 
 \hline
 SGD            & \bf 2.50$\bf \cdot10^{2}$ & 2.30$\cdot10^{2}$ & 2.28$\cdot10^{2}$ & 2.27$\cdot10^{2}$ & 2.26$\cdot10^{2}$ & 2.25$\cdot10^{2}$ & 2.25$\cdot10^{2}$ & 2.24$\cdot10^{2}$ & 2.24$\cdot10^{2}$ & 2.24$\cdot10^{2}$ \\
 \hline
 stoch-smooth     & \bf 2.50$\bf \cdot10^{2}$ & 2.46$\cdot10^{2}$ & 2.34$\cdot10^{2}$ & 2.15$\cdot10^{2}$ & 1.90$\cdot10^{2}$ & 1.61$\cdot10^{2}$ & 1.30$\cdot10^{2}$ & 9.73$\cdot10^{1}$ & 6.77$\cdot10^{1}$ & 4.23$\cdot10^{1}$ \\
 \hline
 Nes-smooth            & \bf 2.50$\bf \cdot10^{2}$ & 1.16$\cdot10^{0}$ & 9.45$\cdot10^{-1}$ & 7.01$\cdot10^{-1}$ & 4.65$\cdot10^{-1}$ & 2.23$\cdot10^{-1}$ & 3.66$\cdot10^{-2}$ & 2.55$\cdot10^{-3}$ & 3.80$\cdot10^{-4}$ & 5.92$\cdot10^{-5}$ \\
 \hline
 TD-$\mu^*$     & \bf 2.50$\bf \cdot10^{2}$ & 1.26$\cdot10^{-2}$ & \bf 2.16$\bf \cdot10^{-9}$ & \bf 2.16$\bf \cdot10^{-9}$ & \bf 2.16$\bf \cdot10^{-9}$ & \bf 2.16$\bf \cdot10^{-9}$ & \bf 2.16$\bf \cdot10^{-9}$ & \bf 2.16$\bf \cdot10^{-9}$ & \bf 2.16$\bf \cdot10^{-9}$ & \bf 2.16$\bf \cdot10^{-9}$ \\
 \hline
 Alg. 1         & \bf 2.50$\bf \cdot10^{2}$ & \bf 2.16$\bf \cdot10^{-9}$ & \bf 2.16$\bf \cdot10^{-9}$ & \bf 2.16$\bf \cdot10^{-9}$ & \bf 2.16$\bf \cdot10^{-9}$ & \bf 2.16$\bf \cdot10^{-9}$ & \bf 2.16$\bf \cdot10^{-9}$ & \bf 2.16$\bf \cdot10^{-9}$ & \bf 2.16$\bf \cdot10^{-9}$ & \bf 2.16$\bf \cdot10^{-9}$ \\
\end{tabular}
\end{adjustbox}
\end{table*}

$\mathcal{R}(y) = \|y\|^2,\, \eta = 0.05$, dimension of matrix $= 300$

\begin{table*}[h!]
\begin{adjustbox}{margin=-40pt 0pt 0pt 0mm} 
\scriptsize
\begin{tabular}{l|l|l|l|l|l|l|l|l|l|l}
 \multicolumn{1}{l|}{Algorithm} & point \#1 & point \#2 & point \#3 & point \#4 & point \#5 & point \#6 & point \#7 & point \#8 & point \#9 & point \#10 \\ 
  \hline
 SGD            & \bf 1.01 & 2.44$\cdot10^{-1}$ & 2.12$\cdot10^{-1}$ & 1.96$\cdot10^{-1}$ & 1.85$\cdot10^{-1}$ & 1.77$\cdot10^{-1}$ & 1.70$\cdot10^{-1}$ & 1.65$\cdot10^{-1}$ & 1.61$\cdot10^{-1}$ & 1.57$\cdot10^{-1}$ \\
  \hline
 stoch-smooth     & \bf 1.01 & 1.01$\cdot10^{0}$ & 1.01$\cdot10^{0}$ & 1.00$\cdot10^{0}$ & 9.92$\cdot10^{-1}$ & 9.79$\cdot10^{-1}$ & 9.65$\cdot10^{-1}$ & 9.47$\cdot10^{-1}$ & 9.28$\cdot10^{-1}$ & 9.06$\cdot10^{-1}$ \\
  \hline
 Nes-smooth             & \bf 1.01 & 2.35$\cdot10^{-1}$ & 9.90$\cdot10^{-3}$ & 7.45$\cdot10^{-4}$ & 3.01$\cdot10^{-3}$ & 1.15$\cdot10^{-3}$ & 5.71$\cdot10^{-6}$ & 3.45$\cdot10^{-4}$ & 3.67$\cdot10^{-4}$ & 4.55$\cdot10^{-5}$ \\
  \hline
 TD-$\mu^*$     & \bf 1.01 & 1.38$\cdot10^{-5}$ & 1.01$\cdot10^{-5}$ & 7.66$\cdot10^{-6}$ & 6.10$\cdot10^{-6}$ & 5.05$\cdot10^{-6}$ & 4.31$\cdot10^{-6}$ & 3.76$\cdot10^{-6}$ & 3.33$\cdot10^{-6}$ & 3.00$\cdot10^{-6}$ \\
  \hline
 Alg. 1         & \bf 1.01 & \bf 1.30$\bf \cdot10^{-7}$ & \bf 3.13$\bf \cdot10^{-8}$ & \bf 1.35$\bf \cdot10^{-8}$ & \bf 7.36$\bf \cdot10^{-9}$ & \bf 4.55$\bf \cdot10^{-9}$ & \bf 3.04$\bf \cdot10^{-9}$ & \bf 2.14$\bf \cdot10^{-9}$ & \bf 1.57$\bf \cdot10^{-9}$ & \bf 1.18$\bf \cdot10^{-9}$ \\
\end{tabular}
\end{adjustbox}
\end{table*}
\bibliographystyle{plain}
\bibliography{references.bib}
\end{document}

%% file: commands.tex
\newtheorem{theorem}{Theorem}[section]

\newtheorem{lemma}[theorem]{Lemma}

\newtheorem{remark}[theorem]{Remark}

%% file: style_PM.tex
\makeatletter
\def\@settitle{\begin{center}%
		\baselineskip14\p@\relax
		\normalfont\LARGE\scshape\bfseries
		\@title
	\end{center}%
}
\makeatother

\makeatletter

\def\subsection{\@startsection{subsection}{2}%
	\z@{.5\linespacing\@plus.7\linespacing}{.5\linespacing}%
	{\normalfont\large\bfseries}}

\def\subsubsection{\@startsection{subsubsection}{3}%
	\z@{.5\linespacing\@plus.7\linespacing}{.5\linespacing}%
	{\normalfont\itshape}}

\allowdisplaybreaks
\date{\today}

\patchcmd\@setauthors
  {\MakeUppercase{\authors}}
  {\authors}
  {}{}